\documentclass[11pt,twoside]{amsart}
\usepackage{graphicx}
\usepackage{color}
\usepackage{amssymb}
\usepackage{amsmath}
\usepackage{amsthm}
\usepackage{amsfonts}
\usepackage{amsmath, amsfonts, amssymb}
\newtheorem{de}{Definition}

\newtheorem{theor}{Theorem}

\numberwithin{equation}{section}

\begin{document}
			\title{Complete classification of two-dimensional associative and diassociative algebras over any basic field}
		\author{I.S.Rakhimov}
		
		\thanks{}
		\maketitle
		\begin{center}
			\address{Universiti Teknologi MARA (UiTM), Shah Alam, Malaysia }
			\end{center}
		
			\begin{abstract}
		A complete classifications, up to isomorphism, of two-dimensional associative and diassociative algebras over any basic field are given.		
	\end{abstract}

	% ----------------------------------------------------------------
	\section{Introduction}
In 1993, Loday introduced the notion of Leibniz algebra \cite{LL}, which
is a generalization of Lie algebra, where the skew-symmetric of the bracket
is dropped and the Jacobi identity is changed by the Leibniz identity. Loday noted that the link between Lie algebras and associative algebras can be extended
to an ``analogous'' link between Leibniz algebra and so-called dialgebra
which is a generalization of associative algebra possessing two products.
Namely, it was shown that if one has a dialgebra $(D, \dashv, \vdash)$ over a finite-dimensional vector space $V$, with two bilinear binary operations with certain compatibility axioms  then introducing a binary operation $[x, y] := x \dashv y-y \vdash x$ we get an algebra structure on $V$ called Leibniz algebra. It also been shown that the universal enveloping algebra of a Leibniz algebra has the structure of a dialgebra.

%In the end of last century Ayupov Sh.A and Omirov B.A \cite{AO}, \cite{OA}, commenced to study of structural properties of Leibniz algebras. They gave the list of isomorphism classes of three dimensional complex Leibniz algebras. Later on in \cite{ARO} one obtained the classification of four dimensional nilpotent non Lie Leibniz algebras.
%
%Leibniz algebras introduced by J-L.Loday in 1993 as a non-antisymmetric version of Lie algebras, whose bracket satisfies the Leibniz identity:
%\begin{center}
%$[[x,y],z]=[[x,z],y]+[x,[y,z]]$
%\end{center}
%
%The Leibniz identity, combined with antisymmetricity, is a variation of the Jacobi identity. Hence Lie algebras are antisymmetric Leibniz algebras.
%Leibniz algebras appeared and investigated in several papers, but to explicitly construct more general examples, we recall that they are related to dialgebras \cite{LF}, in a similar way as Lie algebras are related to associative algebras.\\
%
%A dissociative algebra is a vector space with two bilinear operators $\vdash$, $\dashv$, satisfying five conditions \cite{LF}. Given a dialgebra, the associated Leibniz algebra is obtained by defining the bracket as $[x,y]=x\dashv y - y\vdash x $. Associative algebras are diassociative algebras when the two operations coincide.
In fact, the main motivation of J.-L.Loday to introduce several classes of algebras was the search of an ``obstruction'' to the periodicity in algebraic $K$-theory. 

Since then the study of different properties, relations and classification of Loday's algebras became an active research area. Dozens of papers have been published (see References). But most of the results concerned Loday's algebras over the field of complex numbers. Recently, a result on classification of all algebra structures on two-dimensional vector space over any basic field was published \cite{UB}. In this paper we use the result of \cite{UB} to classify all associative and diassociative algebra structures on two-dimensional vector space over any basic field. This technique was implemented earlier in series of papers \cite{IJAC, MJMS, LJM, UBm} and others. However, there still was a condition on the basic field that was managed to be released in \cite{UB}.

\subsection{Algebras}
\begin{de} A vector space  $\mathbb{V}$ over a field $\mathbb{F}$ equipped with a function
$\cdot :\mathbb{V}\otimes \mathbb{V}\rightarrow \mathbb{V}$
$((\mathrm{x},\mathrm{y})\mapsto \mathrm{x}\cdot \mathrm{y})$ such that
\[(\alpha\mathrm{x}+\beta\mathrm{y})\cdot \mathrm{z}=\alpha(\mathrm{x}\cdot
\mathrm{z})+\beta(\mathrm{y}\cdot \mathrm{z}),\ \ \mathrm{z}\cdot
(\alpha\mathrm{x}+\beta\mathrm{y})=\alpha(\mathrm{z}\cdot
\mathrm{x})+\beta(\mathrm{z}\cdot \mathrm{y})\] whenever $\mathrm{x}, \mathrm{y},
\mathrm{z}\in \mathbb{V}$ and $\alpha, \beta\in \mathbb{F}$, is said to be an
algebra $\mathbb{A}=(\mathbb{V},\cdot)$.\end{de}

\begin{de} \label{ISO} Two algebras $\mathbb{A}$ and $\mathbb{B}$ are called isomorphic
if there is an invertible linear map  $f:\mathbb{A}\rightarrow \mathbb{B} $
such that \[f(\mathrm{x}\cdot_{\mathbb{A}}
\mathrm{y})=f(\mathrm{x})\cdot_{\mathbb{B}} f(\mathrm{y})\]
whenever $\mathrm{x}, \mathrm{y}\in \mathbb{A}.$\end{de}

\begin{de} An invertible linear map  $f:\mathbb{A}\rightarrow
\mathbb{A} $ is said to be an automorphism if \[f(\mathrm{x}\cdot
\mathrm{y})=f(\mathrm{x})\cdot f(\mathrm{y})\] whenever
$\mathrm{x}, \mathrm{y}\in \mathbb{A}.$ \end{de}

The set of all automorphisms of an algebra $\mathbb{A}$ forms a group with respect to the composition operation and it is denoted by $Aut(\mathbb{A}).$

Let $\mathbb{A}$ be $n$-dimensional algebra over $\mathbb{F}$ and
$\mathbf{e}=(\mathrm{e}_1,\mathrm{e}_2,...,\mathrm{e}_n)$ be its basis. Then the bilinear map $\cdot$ is represented by a $n \times n^2$ matrix (called the matrix of structure constant, shortly MSC) $$A=\left(\begin{array}{ccccccccccccc}a_{11}^1&a_{12}^1&...&a_{1n}^1&a_{21}^1&a_{22}^1&...&a_{2n}^1&...&a_{n1}^1&a_{n2}^1&...&a_{nn}^1\\ a_{11}^2&a_{12}^2&...&a_{1n}^2&a_{21}^2&a_{22}^2&...&a_{2n}^2&...&a_{n1}^2&a_{n2}^2&...&a_{nn}^2 \\
...&...&...&...&...&...&...&...&...&...&...&...&...\\ a_{11}^n&a_{12}^n&...&a_{1n}^n&a_{21}^n&a_{22}^n&...&a_{2n}^n&...&a_{n1}^n&a_{n2}^n&...&a_{nn}^n\end{array}\right)$$ as follows
$$\mathrm{e}_i \cdot \mathrm{e}_j=\sum\limits_{k=1}^n a_{ij}^k\mathrm{e}_k, \ \mbox{where}\ i,j=1,2,...,n.$$
Therefore, the product on $\mathbb{A}$ with respect to the basis $\mathbf{e}$ is written as follows
\begin{equation} \label{Product}
	\mathrm{x}\cdot \mathrm{y}=\mathbf{e} A(x\otimes y)
\end{equation}
 for any $\mathrm{x}=\mathbf{e}x,\mathrm{y}=\mathbf{e}y,$
	where $x=(x_1, x_2,...,x_n)^T,$ and  $y=(y_1, y_2,...,y_n)^T$ are column coordinate vectors of $\mathrm{x}$ and $\mathrm{y},$ respectively, $x\otimes y$ is the tensor(Kronecker)
 product of the vectors $x$ and $y$. Now and onward for the product ``$\mathrm{x}\cdot \mathrm{y}$'' on $\mathbb{A}$ we use the juxtaposition ``$\mathrm{x} \mathrm{y}$''.

 Further we
assume that the basis $\mathbf{e}$ is fixed and we do not make a difference between the algebra
$\mathbb{A}$ and its MSC $A$.

An automorphism $\mathrm{g}:\mathbb{A}\rightarrow \mathbb{A}$ as an invertible linear
map is represented on the basis $\mathbf{e}$ by an invertible matrix $g\in GL(n;\mathbb{F})$ and
$\mathrm{g}(\mathrm{x})=\mathrm{g}(\mathbf{e}x)=\mathbf{e}gx$. Due to  \[\mathrm{g}(\mathrm{x}\cdot
\mathrm{y})=\mathrm{g}(\mathbf{e}A(x\otimes y))=\mathbf{e}g(A(x\otimes y))=\mathbf{e}(gA)(x\otimes y),\] and
\[\mathrm{g}(\mathrm{x})\cdot \mathrm{g}(\mathrm{y})=(\mathbf{e}gx)\cdot (\mathbf{e}gy)=\mathbf{e}A(gx\otimes
gy)=\mathbf{e}Ag^{\otimes 2}(x\otimes y)\] the condition $\mathrm{g}(\mathrm{x}\cdot
\mathrm{y})=\mathrm{g}(\mathrm{x})\cdot \mathrm{g}(\mathrm{y})$ is written in terms of $A$ and $g$ as follows
\begin{equation}\label{1}gA=Ag^{\otimes 2}.\end{equation}

Note that in this term Definition \ref{ISO} can also be rewritten as
\begin{equation}\label{ISO1}gA=Bg^{\otimes 2} \Longleftrightarrow A=g^{-1}Bg^{\otimes 2}.\end{equation}

%An derivation $\mathrm{d}:\mathbb{A}\rightarrow \mathbb{A}$ as a linear map is
%represented by a matrix $d\in M(m;\mathbb{F})$ as follows
%$\mathrm{d}(\mathrm{x})=\mathrm{d}(eu)=edu$. Due to  \[\mathrm{d}(\mathrm{x}\cdot
%\mathrm{y})=\mathrm{d}(eA(u\otimes v))=ed(A(u\otimes v))=e(dA)(u\otimes v)\] and
%\[\mathrm{d}(\mathrm{x})\cdot \mathrm{y}+\mathrm{x}\cdot
%\mathrm{d}(\mathrm{y})=(edu)\cdot (ev)+(eu)\cdot (edv)=eA(du\otimes v)+eA(u\otimes
%dv)=\]
%\[e(A(d\otimes I)(u\otimes v)+A(I\otimes d)(u\otimes v))= eA((d\otimes I)+(I\otimes
%d))(u\otimes v)\]  the property $\mathrm{d}(\mathrm{x}\cdot
%\mathrm{y})=\mathrm{d}(\mathrm{x})\cdot \mathrm{y}+\mathrm{x}\cdot
%\mathrm{d}(\mathrm{y})$ is equivalent to \begin{equation}\label{2}dA=A(d\otimes
%I+I\otimes d),\end{equation}
%where $I$ stands for the identity matrix.
\subsection{Associative algebras}	
 \begin{de}  An algebra $(\mathbb{A},\cdot)$ is said to be associative if
 %\begin{center}
%$\succ:E \times E \rightarrow E$,\; $\prec:E \times E \rightarrow E$
%\end{center}
$\forall\ \mathrm{x}, \mathrm{y}, \mathrm{z} \in \mathbb{A}$ the following axiom holds true
\begin{equation} \label{AA}
\begin{array}{rll}
(\mathrm x \cdot \mathrm y)\cdot \mathrm z &=& \mathrm x \cdot (\mathrm y \cdot \mathrm z),\\
%\mathrm x \dashv (\mathrm y \dashv \mathrm z) &=&  \mathrm x \dashv(\mathrm y \vdash \mathrm z),\\
%(\mathrm x \vdash \mathrm y)\dashv \mathrm z &=& \mathrm x \vdash(\mathrm y \dashv \mathrm z),\\
%(\mathrm x \dashv \mathrm y)\vdash \mathrm z &=& (\mathrm x \vdash \mathrm y)\vdash \mathrm z,\\
%(\mathrm x \vdash \mathrm y)\vdash \mathrm z &=& \mathrm x \vdash(\mathrm y \vdash \mathrm z).\\
\end{array}
\end{equation}
 \end{de}

Write  $$\mathrm{x}\cdot \mathrm{y}=\mathbf{e}A(x\otimes y) \ \mbox{and} \ \mathrm{y}\cdot \mathrm{z}=\mathbf{e}A(y\otimes z),$$
$$(\mathrm{x}\cdot \mathrm{y})\cdot \mathrm{z}=\mathbf{e}A(A(x\otimes y)\otimes z))\ \mbox{and} \ \mathrm{x}\cdot (\mathrm{y}\cdot \mathrm{z})=\mathbf{e}A(x\otimes A(y\otimes z)).$$
	%where $\mathbf{e}=\{\mathrm{e}_1,\mathrm{e}_2,...,\mathrm{e}_n\}$ is a basis of $(\mathbb{A}$, $A$ is matrix of the structure constants of $\mathbb{A}$ with respect to the binary operation $\cdot$, $x=(x_1, x_2,...,x_2)^T,$ $y=(y_1, y_2,...,y_n)^T$ and $z=(z_1, z_2,...,z_n)^T$ are column coordinate vectors of $\mathrm{x}$, $\mathrm{y}$  and $\mathrm{z}$, respectively, and $\otimes $ is the tensor (Kronecker)
% product.
	
	Then,

%$$\mathrm x \dashv(\mathrm y \vdash \mathrm z)=\mathbf{e}A(x\otimes B(y\otimes z)),$$
%$$(\mathrm{x}\vdash \mathrm{y})\dashv \mathrm{z}=\mathbf{e}A(B(x\otimes y)\otimes z)),$$
%$$\mathrm{x}\vdash (\mathrm{y}\dashv \mathrm{z})=\mathbf{e}B(x\otimes (A(y\otimes z)),$$
%$$(\mathrm{x}\dashv \mathrm{y})\vdash \mathrm{z}=\mathbf{e}B(A(x\otimes y)\otimes z),$$
%$$(\mathrm{x}\vdash \mathrm{y})\vdash \mathrm{z}=\mathbf{e}B(B(x\otimes y)\otimes z)),$$
%$$\mathrm{x}\vdash (\mathrm{y}\vdash \mathrm{z})=\mathbf{e}B(x\otimes (B(y\otimes z)),$$

$$\mathbf{e}A(A(x\otimes y)\otimes z))=\mathbf{e}A(x\otimes A(y\otimes z)),$$
i.e., an algebra $\mathbb{A}$ with MSC $A$ is associative if and only if
\begin{equation} \label{AA2}
\begin{array}{ccc}
A(A\otimes I)&=&A(I\otimes A)
\end{array},
\end{equation}
where $I$ is $n\times n$ identity matrix.
\subsection{Associative dialgebras}	
\begin{de}  A dialgebra $\mathbb{D}=(\mathbb{V},\dashv,\vdash)$ is said to be associative dialgebra if
 %\begin{center}
%$\succ:E \times E \rightarrow E$,\; $\prec:E \times E \rightarrow E$
%\end{center}
$\forall \mathrm{x}, \mathrm{y}, \mathrm{z} \in \mathbb{D}$ the following axioms hold true
\begin{equation} \label{DAA}
\begin{array}{rll}
(\mathrm x \dashv \mathrm y)\dashv \mathrm z &=& \mathrm x \dashv (\mathrm y \dashv \mathrm z),\\
\mathrm x \dashv (\mathrm y \dashv \mathrm z) &=&  \mathrm x \dashv(\mathrm y \vdash \mathrm z),\\
(\mathrm x \vdash \mathrm y)\dashv \mathrm z &=& \mathrm x \vdash(\mathrm y \dashv \mathrm z),\\
(\mathrm x \dashv \mathrm y)\vdash \mathrm z &=& (\mathrm x \vdash \mathrm y)\vdash \mathrm z,\\
(\mathrm x \vdash \mathrm y)\vdash \mathrm z &=& \mathrm x \vdash(\mathrm y \vdash \mathrm z).\\
\end{array}
\end{equation}
 \end{de}

 \begin{de} Let $\mathbb{D}_1=(\mathbb{V},\dashv,\vdash)$ and $\mathbb{D}_2=(\mathbb{V},\dashv',\vdash')$ be diassociative algebras. A linear function $f:\mathbb{D}_1\longrightarrow \mathbb{D}_2$ is said to be homomorphism if
 $$f(\mathrm x \dashv \mathrm y)=f(\mathrm x) \dashv' f(\mathrm y) \ \mbox{and}\  f(\mathrm x \vdash \mathrm y)=f(\mathrm x) \vdash' f(\mathrm y) \ \mbox{for all}\ \mathrm x, \mathrm y\in \mathbb{D}_1. $$
 \end{de}
  \begin{de}
 Dialgebras  $\mathbb{D}_1$ and $\mathbb{D}_2$ are called isomorphic
if there is an invertible homomorphism  $f:\mathbb{D}_1\longrightarrow \mathbb{D}_2$.
\end{de}

Let  \[\mathrm{x}\dashv \mathrm{y}=\mathbf{e}A(x\otimes y)\ \mbox{and} \ \mathrm{x}\vdash \mathrm{y}=\mathbf{e}B(x\otimes y)\] for any $\mathrm{x}=\mathbf{e}x,\mathrm{y}=\mathbf{e}y.$
	%where $\mathbf{e}=\{\mathrm{e}_1,\mathrm{e}_2\}$ is a basis of $\mathbf{D}$, $A$ and $B$ are matrices of the structure constants of $\mathbf{D}$ with respect to the operations $\dashv$ and $\vdash$, respectively, $x=(x_1, x_2)^T,$ and  $y=(y_1, y_2)^T$ are column coordinate vectors of $\mathrm{x}$ and $\mathrm{y},$ respectively, $(x\otimes y)=(x_1y_1,x_1y_2,x_2y_1,x_2y_2)^T$ is the tensor (Kronecker)
% product of the vectors $x$ and $y$.
	
	Then,
$$(\mathrm{x}\dashv \mathrm{y})\dashv \mathrm{z}=\mathbf{e}A(A(x\otimes y)\otimes z)),$$
%\hspace{3cm}
$$\mathrm{x}\dashv (\mathrm{y}\dashv \mathrm{z})=\mathbf{e}A(x\otimes A(y\otimes z)),$$
$$\mathrm x \dashv(\mathrm y \vdash \mathrm z)=\mathbf{e}A(x\otimes B(y\otimes z)),$$
$$(\mathrm{x}\vdash \mathrm{y})\dashv \mathrm{z}=\mathbf{e}A(B(x\otimes y)\otimes z)),$$
$$\mathrm{x}\vdash (\mathrm{y}\dashv \mathrm{z})=\mathbf{e}B(x\otimes (A(y\otimes z)),$$
$$(\mathrm{x}\dashv \mathrm{y})\vdash \mathrm{z}=\mathbf{e}B(A(x\otimes y)\otimes z),$$
$$(\mathrm{x}\vdash \mathrm{y})\vdash \mathrm{z}=\mathbf{e}B(B(x\otimes y)\otimes z)),$$
$$\mathrm{x}\vdash (\mathrm{y}\vdash \mathrm{z})=\mathbf{e}B(x\otimes (B(y\otimes z)),$$

%and
%$$\mathrm x \dashv (\mathrm y \dashv \mathrm z) + \mathrm x \dashv(\mathrm y \vdash \mathrm z)=\mathbf{e}(A(x\otimes A(y\otimes z))+A(x\otimes B(y\otimes z))),$$
%$$(\mathrm x \dashv \mathrm y)\vdash \mathrm z + (\mathrm x \vdash \mathrm y)\vdash \mathrm z=\mathbf{e}(B(A(x\otimes y)\otimes z)+B(B(x\otimes y)\otimes z))).$$
 %\mathbf{y}^2=eAy^{\otimes 2}, \mathbf{x}^2\mathbf{y}^2=eA(Ax^{\otimes 2}\otimes Ay^{\otimes 2})=eA(A\otimes A)(x^{\otimes 2}\otimes y^{\otimes 2})$$ and
%	$$(\mathbf{x}\cdot \mathbf{y})^2=eA(A(x\otimes y)\otimes A(x\otimes y))=eA(A\otimes A)(x\otimes y)^{\otimes 2}.$$
%	Therefore, the equality $\mathbf{x}^2\mathbf{y}^2=(\mathbf{x}\cdot \mathbf{y})^2$ can be written in terms of MSC as follows
%	$$eA(A\otimes A)(x^{\otimes 2}\otimes y^{\otimes 2})=eA(A\otimes A)(x\otimes y)^{\otimes 2}$$ and
%	$$A(A\otimes A)(x^{\otimes 2}\otimes y^{\otimes 2}-(x\otimes y)^{\otimes 2})=0$$ or just \[A(A\otimes A)(x\otimes x\otimes y\otimes y-x\otimes y \otimes x\otimes y)=0.\ \ \ \ \ *\]
Therefore, the diassociative algebra axioms (\ref{DAA}) in terms of the structure constants can be given by the identities
\begin{equation} \label{DAA2}
\begin{array}{rll}
 A(A(x\otimes y)\otimes z)&=&A(x\otimes A(y\otimes z)),\\
 A(x\otimes A(y\otimes z))&=&A(x\otimes B(y\otimes z)),\\
A(B(x\otimes y)\otimes z)&=&B(x\otimes (A(y\otimes z)),\\
B(A(x\otimes y)\otimes z)&=&B(B(x\otimes y)\otimes z)),\\
B(B(x\otimes y)\otimes z))&=&B(x\otimes (B(y\otimes z)).
\end{array}
\end{equation}	

The axioms can be rewritten as follows
\begin{equation} \label{DAA3}
\begin{array}{rll}
 A(A \otimes I)&-&A(I \otimes A)=0\\
 A(I \otimes A)&-&A(I \otimes B)=0\\
A(B\otimes I)&-&B(I \otimes A)=0\\
B(A \otimes I)&-&B(B \otimes I)=0\\
B(B \otimes I)&-&B  (B \otimes I)=0
\end{array}
\end{equation}	
i.e., a dialgebra $\mathbb{D}$ with MSC $D:=\{A, B\}$ is diassociative if and only if (\ref{DAA3}) holds true.

In the paper we make use the following result from \cite{UB} on complete classification of two-dimensional algebras over any basic field.

\begin{theor} \label{Char0} Any non-trivial $2$-dimensional algebra over a field $\mathbb{F}$ $(Char(\mathbb{F})\neq 2,3)$ is isomorphic to only one of the following listed, by their matrices of structure constants, such algebras:
 		\begin{itemize}
 			\item $A_1(\mathrm{c})=\begin{pmatrix}
 			\alpha_1&\alpha_2&1+\alpha_2&\alpha_4 \\
 			\beta_1& -\alpha_1& 1-\alpha_1& -\alpha_2
 			\end{pmatrix},$
 			where $ \mathrm{c}=(\alpha_1,\alpha_2,\alpha_4,\beta_1)\in\mathbb{F}^4,$
 			
 			\item $A_2(\mathrm{c})=\begin{pmatrix}
 			\alpha_1&0&0&\alpha_4 \\
 			1& \beta_2& 1-\alpha_1& 0
 			\end{pmatrix},$
 			where $ \mathrm{c}=(\alpha_1,\alpha_4,\beta_2)\in\mathbb{F}^3$ and $\alpha_4\neq0,$
 			
 			\item $A_3(\mathrm{c})=\begin{pmatrix}
 			\alpha_1&0&0&\alpha_4 \\
 			0& \beta_2& 1-\alpha_1& 0
 			\end{pmatrix}\simeq\begin{pmatrix}
 			\alpha_1&0&0&a^2\alpha_4 \\
 			0& \beta_2& 1-\alpha_1& 0
 			\end{pmatrix},$  where $ \mathrm{c}=(\alpha_1,\alpha_4,\beta_2)\in\mathbb{F}^3,$ $a\in\mathbb{F}$ and $ a\neq 0,$
 			
 			\item $A_4(\mathrm{c})=\begin{pmatrix}
 			0&1&1&0 \\
 			\beta_1& \beta_2& 1& -1
 			\end{pmatrix},$  where $ \mathrm{c}=(\beta_1,\beta_2)\in\mathbb{F}^2,$
 			
 			\item $A_5(\mathrm{c})=\begin{pmatrix}
 			\alpha_1&0&0&0 \\
 			1& 2\alpha_1-1& 1-\alpha_1& 0
 			\end{pmatrix},$\\ where $ \mathrm{c}=\alpha_1\in\mathbb{F},$
 			
 			\item $A_6(\mathrm{c})=\begin{pmatrix}
 			\alpha_1&0&0&\alpha_4 \\
 			1& 1-\alpha_1& -\alpha_1& 0
 			\end{pmatrix},$ where $ \mathrm{c}=(\alpha_1,\alpha_4)\in\mathbb{F}^2$ and $ \alpha_4\neq0,$
 			
 			\item $A_7(\mathrm{c})=\begin{pmatrix}
 			\alpha_1&0&0&\alpha_4 \\
 			0& 1-\alpha_1& -\alpha_1& 0
 			\end{pmatrix}\simeq\begin{pmatrix}
 			\alpha_1&0&0&a^2\alpha_4 \\
 			0& 1-\alpha_1& -\alpha_1& 0
 			\end{pmatrix},$ where $ \mathrm{c}=(\alpha_1,\alpha_4)\in\mathbb{F}^2, $ $ a\in\mathbb{F}$ and $ a\neq 0,$
 			
 			\item $A_8(\mathrm{c})=\begin{pmatrix}
 			0&1&1&0 \\
 			\beta_1& 1&0& -1
 			\end{pmatrix},$ where $ \mathrm{c}=\beta_1\in\mathbb{F},$
 			
 			\item $A_9=\begin{pmatrix}
 			\frac{1}{3}&0&0&0\\
 			1&\frac{2}{3}&-\frac{1}{3}&0
 			\end{pmatrix},
 			$
 			
 			\item $A_{10}(\mathrm{c})=\begin{pmatrix}
 			0&1&1&1\\
 			\beta_1&0&0&-1
 			\end{pmatrix}\simeq\begin{pmatrix}
 			0&1&1&1\\
 			\beta_1^{'}(a)&0&0&-1
 			\end{pmatrix},$\\ where $ \mathrm{c}=\beta_1\in\mathbb{F},$ the polynomial $(\beta_1t^3-3t-1)(\beta_1t^2+\beta_1t+1)(\beta_1^2t^3+6\beta_1t^2+3\beta_1t+\beta_1-2)$ has no root in $\mathbb{F}$,
 			$a\in\mathbb{F}$ and $ \beta_1^{'}(t)=\frac{(\beta_1^2t^3+6\beta_1t^2+3\beta_1t+\beta_1-2)^2}{(\beta_1t^2+\beta_1t+1)^3},$
 			
 			\item $A_{11}(\mathrm{c})=\begin{pmatrix}
 			0&0&0&1 \\
 			\beta_1&0&0&0\\
 			\end{pmatrix}\simeq\begin{pmatrix}
 			0&0&0&1 \\
 			a^{3}\beta_1^{\pm1}&0&0&0\\
 			\end{pmatrix},$ where the polynomial $ \beta_1-t^3$ has no root in $\mathbb{F},$\ $ a, \mathrm{c}=\beta_1\in \mathbb{F}$  and $ a, \beta_1\neq0,$
 			
 			\item $A_{12}(\mathrm{c})=\begin{pmatrix}
 			0&1&1&0 \\
 			\beta_1&0&0&-1
 			\end{pmatrix}\simeq\begin{pmatrix}
 			0&1&1&0 \\
 			a^2\beta_1&0&0&-1
 			\end{pmatrix},$
 			where $a, \mathrm{c}=\beta_1\in\mathbb{F}$ and $ a\neq 0,$
 			
 			\item $A_{13}=\begin{pmatrix}
 			0&0&0&0 \\
 			1&0&0&0
 			\end{pmatrix}.
 			$
 		\end{itemize}
 	\end{theor}
 \begin{theor} \label{Char2} Any non-trivial $2$-dimensional algebra over a field $\mathbb{F}$ $(Char(\mathbb{F})=2)$ is isomorphic to only one of the following listed by their matrices of structure constants, such algebras:
 	\begin{itemize}
 		\item
 		$ A_{1,2}(\mathrm{c})=\begin{pmatrix}
 		\alpha_1&\alpha_2&1+\alpha_2&\alpha_4 \\
 		\beta_1& \alpha_1& 1+\alpha_1&\alpha_2
 		\end{pmatrix},
 		$ where $ \mathrm{c}=(\alpha_1,\alpha_2,\alpha_4,\beta_1)\in\mathbb{F}^4$
 		
 		\item
 		$ A_{2,2}(\mathrm{c})=\begin{pmatrix}
 		\alpha_1&0&0&\alpha_4 \\
 		1& \beta_2& 1+\alpha_1& 0
 		\end{pmatrix},
 		$ where $ \mathrm{c}=(\alpha_1,\alpha_4,\beta_2)\in\mathbb{F}^3$ and $\alpha_4\neq 0$
 		\item
 		$ A_{2,2}(\alpha_1,0,1)=\begin{pmatrix}
 		\alpha_1&0&0&0 \\
 		1& 1& 1+\alpha_1& 0
 		\end{pmatrix},
 		$ where $\alpha_1\in\mathbb{F}$
 		
 		\item $ A_{3,2}(\mathrm{c})=\begin{pmatrix}
 		\alpha_1&0&0&\alpha_4 \\
 		0& \beta_2& 1+\alpha_1& 0
 		\end{pmatrix}\simeq\begin{pmatrix}
 		\alpha_1&0&0&a^2\alpha_4 \\
 		0& \beta_2& 1+\alpha_1& 0
 		\end{pmatrix}
 		,$ where $ \mathrm{c}=(\alpha_1,\alpha_4,\beta_2)\in\mathbb{F}^3,$ $a\in\mathbb{F}$ and $ a\neq 0$

 		\item $ A_{4,2}(\mathrm{c})=\begin{pmatrix}
 		\alpha_1&1&1&0 \\
 		\beta_1& \beta_2& 1+\alpha_1& 1
 		\end{pmatrix}\simeq\begin{pmatrix}
 		\alpha_1&1&1&0 \\
 		\beta_1+(1+\beta_2)a+a^2& \beta_2& 1+\alpha_1& 1
 		\end{pmatrix},
 		$\\ where $ \mathrm{c}=(\alpha_1,\beta_1,\beta_2)\in\mathbb{F}^3$
 		
 		\item $A_{5,2}(\mathrm{c})=\begin{pmatrix}
 		\alpha_1&0&0&\alpha_4 \\
 		1& 1+\alpha_1& \alpha_1& 0
 		\end{pmatrix},
 		$ where $ \mathrm{c}=(\alpha_1,\alpha_4)\in\mathbb{F}^2$ and $\alpha_4\neq0$
 		\item $A_{5,2}(1,0)=\begin{pmatrix}
 		1&0&0&0 \\
 		1&0&1& 0
 		\end{pmatrix},
 		$
 		
 		\item $A_{6,2}(\mathrm{c})=\begin{pmatrix}
 		\alpha_1&0&0&\alpha_4 \\
 		0& 1+\alpha_1& \alpha_1& 0
 		\end{pmatrix}\simeq\begin{pmatrix}
 		\alpha_1&0&0&a^2\alpha_4 \\
 		0& 1+\alpha_1& \alpha_1& 0
 		\end{pmatrix}
 		,$ where $ \mathrm{c}=(\alpha_1,\alpha_4)\in\mathbb{F}^2,$ $a\in \mathbb{F}$ and $a\neq0$
 		
 		\item $A_{7,2}(\mathrm{c})=\begin{pmatrix}
 		\alpha_1&1&1&0 \\
 		\beta_1& 1+\alpha_1& \alpha_1& 1
 		\end{pmatrix}\simeq\begin{pmatrix}
 		\alpha_1&1&1&0 \\
 		\beta_1+a\alpha_1+a+a^2& 1+\alpha_1& \alpha_1& 1
 		\end{pmatrix}
 		,$ where $ \mathrm{c}=(\alpha_1,\beta_1)\in\mathbb{F}^2$ and $a\in\mathbb{F}$
 		
 		\item $A_{8,2}(\mathrm{c})=\begin{pmatrix}
 		0&1&1&1 \\
 		\beta_1& 0&0& 1
 		\end{pmatrix}\simeq\begin{pmatrix}
 		0&1&1&1 \\
 		\beta_1^{'}(a)& 0&0& 1
 		\end{pmatrix}
 		,$ where the polynomial\\ $(\beta_1t^3+t+1)(\beta_1t^2+\beta_1t+1)$ has no root in $\mathbb{F},$
 		$a\in\mathbb{F}$ and $ \beta_1^{'}(t)=\frac{(\beta_1^2t^3+\beta_1t+\beta_1)^2}{(\beta_1t^2+\beta_1t+1)^3}$
 		
 		\item $A_{9,2}(\mathrm{c})=\begin{pmatrix}
 		0&0&0&1\\
 		\beta_1&0&0&0
 		\end{pmatrix}\simeq\begin{pmatrix}
 		0&0&0&1\\
 		a^3\beta_1^{\pm1}&0&0&0
 		\end{pmatrix}
 		,$ where $ a, \mathrm{c}=\beta_1\in\mathbb{F}$ and $a\neq 0,$
 		
 		the polynomial $ \beta_1+t^3$ has no root in $\mathbb{F}$
 		
 		\item $A_{10,2}(\mathrm{c})=\begin{pmatrix}
 		1&1&1&0\\
 		\beta_1&1&1&1
 		\end{pmatrix}\simeq\begin{pmatrix}
 		1&1&1&0\\
 		\beta_1+a+a^2&1&1&1
 		\end{pmatrix}
 		,$ where $a, \mathrm{c}=\beta_1\in\mathbb{F}$
 		
 		\item $A_{11,2}(\mathrm{c})=\begin{pmatrix}
 		0&1&1&0 \\
 		\beta_1&0&0&1\\
 		\end{pmatrix}\simeq\begin{pmatrix}
 		0&1&1&0 \\
 		b^2(\beta_1+a^2)&0&0&1\\
 		\end{pmatrix}
 		,$ where $a, b\in\mathbb{F}$ and $ b\neq 0$
 		
 		\item $A_{12,2}=\begin{pmatrix}
 		0&0&0&0 \\
 		1&0&0&0\\
 		\end{pmatrix}
 		$
  		
 	\end{itemize}
 \end{theor}

 \begin{theor} \label{Char3} Any non-trivial $2$-dimensional algebra over a field $\mathbb{F}$ $(Char(\mathbb{F})=3)$ is isomorphic to only one of the following, listed by their matrices of structure constants, such algebras:
 	\begin{itemize}
 		\item $A_{1,3}(\mathrm{c})=\begin{pmatrix}
 		\alpha_1&\alpha_2&\alpha_2+1&\alpha_4 \\
 		\beta_1& -\alpha_1& 1-\alpha_1& -\alpha_2
 		\end{pmatrix},$ where $ \mathrm{c}=(\alpha_1,\alpha_2,\alpha_4,\beta_1)\in\mathbb{F}^4$
 		
 		\item $A_{2,3}(\mathrm{c})=\begin{pmatrix}
 		\alpha_1&0&0&\alpha_4 \\
 		1& \beta_2& 1-\alpha_1& 0
 		\end{pmatrix},$ where $ \mathrm{c}=(\alpha_1,\alpha_4,\beta_2)\in\mathbb{F}^3,$ and $\alpha_4\neq0$
 		
 		\item $A_{3,3}(\mathrm{c})=\begin{pmatrix}
 		\alpha_1&0&0&\alpha_4 \\
 		0& \beta_2& 1-\alpha_1& 0
 		\end{pmatrix}\simeq\begin{pmatrix}
 		\alpha_1&0&0&a^2\alpha_4 \\
 		0& \beta_2& 1-\alpha_1& 0
 		\end{pmatrix},$ where $ \mathrm{c}=(\alpha_1,\alpha_4,\beta_2)\in\mathbb{F}^3,$ $ a\in\mathbb{F}$ and $ a\neq 0$
 		
 		\item $A_{4,3}(\mathrm{c})=\begin{pmatrix}
 		0&1&1&0 \\
 		\beta_1& \beta_2& 1& -1
 		\end{pmatrix},$ where $ \mathrm{c}=(\beta_1,\beta_2)\in\mathbb{F}^2$
 		
 		\item $A_{5,3}(\mathrm{c})=\begin{pmatrix}
 		\alpha_1&0&0&0 \\
 		1& 2\alpha_1-1& 1-\alpha_1& 0
 		\end{pmatrix},$ where $ \mathrm{c}=\alpha_1\in\mathbb{F}$
 		
 		\item $A_{6,3}(\mathrm{c})=\begin{pmatrix}
 		\alpha_1&0&0&\alpha_4 \\
 		1& 1-\alpha_1& -\alpha_1& 0
 		\end{pmatrix},$ where $ \mathrm{c}=(\alpha_1,\alpha_4)\in\mathbb{F}^2$ and $ \alpha_4\neq0$
 		
 		\item $A_{7,3}(\mathrm{c})=\begin{pmatrix}
 		\alpha_1&0&0&\alpha_4 \\
 		0& 1-\alpha_1& -\alpha_1& 0
 		\end{pmatrix}\simeq\begin{pmatrix}
 		\alpha_1&0&0&a^2\alpha_4 \\
 		0& 1-\alpha_1& -\alpha_1& 0
 		\end{pmatrix},$ where $ \mathrm{c}=(\alpha_1,\alpha_4)\in\mathbb{F}^2,$ $a\in\mathbb{F}$ and $ a\neq 0$
 		
 		\item $A_{8,3}(\mathrm{c})=\begin{pmatrix}
 		0&1&1&0 \\
 		\beta_1& 1&0& -1
 		\end{pmatrix},$ where $ \mathrm{c}=\beta_1\in\mathbb{F}$
 		
 		\item $A_{9,3}(\beta_1)=\begin{pmatrix}
 		0&1&1&1\\
 		\beta_1&0&0&-1
 		\end{pmatrix}\simeq\begin{pmatrix}
 		0&1&1&1\\
 		\beta_1'(a)&0&0&-1
 		\end{pmatrix},$ where the polynomial \\ $(\beta_1-t^3)(\beta_1t^2+\beta_1t+1)(\beta_1^2t^3+\beta_1-2)$ has no root in $\mathbb{F},$
 		$a\in\mathbb{F}$ and $ \beta_1^{'}(t)=\frac{(\beta_1^2t^3+\beta_1-2)^2}{(\beta_1t^2+\beta_1t+1)^3}$
 		
 		\item $A_{10,3}(\mathrm{c})=\begin{pmatrix}
 		0&0&0&1\\
 		\beta_1&0&0&0
 		\end{pmatrix}\simeq\begin{pmatrix}
 		0&0&0&1\\
 		a^{3}\beta_1^{\pm1}&0&0&0
 		\end{pmatrix},$\\  where the polynomial $ \beta_1-t^3$ has no root, $ a, \mathrm{c}=\beta_1\in\mathbb{F}$ and  $ a, \beta_1\neq 0$
 		
 		\item $A_{11,3}(\mathrm{c})=\begin{pmatrix}
 		0&1&1&0 \\
 		\beta_1&0&0&-1\\
 		\end{pmatrix}\simeq\begin{pmatrix}
 		0&1&1&0 \\
 		a^{2}\beta_1&0&0&-1\\
 		\end{pmatrix},
 		$ where  $a, \mathrm{c}=\beta_1\in\mathbb{F},$ $a\neq 0$
 		
 		\item $A_{12,3}=\begin{pmatrix}
 		1&0&0&0 \\
 		1&-1&-1&0
 		\end{pmatrix}
 		,$
 		
 		\item $A_{13,3}=\begin{pmatrix}
 		0&0&0&0 \\
 		1&0&0&0
 		\end{pmatrix}.$
 	\end{itemize}
 \end{theor}

 The next sections are devoted to the classification of all two-dimensional associative and associative dialgebras over any basic field relying on the theorems above.

\section{Classification of two-dimensional associative algebras}
In this section we classify all two-dimensional associative algebras over any basic field.
Let $\mathbb{A}$ be a two-dimensional associative algebra and \[A=\left(\begin{array}{cccc} \alpha_1 & \alpha_2 & \alpha_3 &\alpha_4\\ \beta_1 & \beta_2 & \beta_3 &\beta_4\end{array}\right)\] be its MSC on a basis $\mathbf{e}=\left(\mathrm{e}_1, \mathrm{e}_2\right)$. Write the axiom (\ref{AA2}) in terms of the elements of $A$ as follows

\begin{equation} \label{GSEAs}
\begin{array}{lll}
  \beta_1(\alpha_2-\alpha_3)&=&0\\
  \alpha_2\beta_2-\alpha_4\beta_1&=&0\\
 (\alpha_1-\beta_3)\alpha_2-\alpha_3(\alpha_1-\beta_2)&=&0\\
 (\alpha_1-\beta_2)\alpha_4-\alpha_2(\alpha_2-\beta_4)&=&0\\
  \alpha_3\beta_3-\alpha_4\beta_1&=&0\\
  \alpha_4(\beta_2-\beta_3)&=&0\\
 (\alpha_1-\beta_3)\alpha_4-\alpha_3(\alpha_3-\beta_4)&=&0\\
 \alpha_4(\alpha_2-\alpha_3)&=&0\\
 \beta_1(\beta_2-\beta_3)&=&0\\
 (\alpha_2-\beta_4)\beta_1-\beta_2(\alpha_1-\beta_2)&=&0\\
% \beta_1(\alpha_2-\alpha_3)&=&0\\
% \alpha_2\beta_2-\alpha_4\beta_1&=&0\\
 (\alpha_3-\beta_4)\beta_1-\beta_3(\alpha_1-\beta_3)&=&0\\
 (\alpha_3-\beta_4)\beta_2-\beta_3(\alpha_2-\beta_4)&=&0\\
% \alpha_3\beta_3-\alpha_4\beta_1&=&0\\
% \alpha_4(\beta_2-\beta_3)&=&0\\
%\begin{equation} \label{GSE}
%\begin{array}{lll}
\end{array}
\end{equation}

Theorems \ref{Char0}, \ref{Char2}, \ref{Char3} are applied as follows: substitute the structure constants of the list of representatives in the theorems into the system of equations (\ref{GSEAs}) taking the structure constants to be variables. The solutions to the system give structure constants of associative algebras.

\subsection{Characteristic is not 2 and 3}\

It is easy to see that $A_{13}$ is associative.

For algebras $A_{12} - A_4$ the system of equations (\ref{GSEAs}) is inconsistent.

Consider

$A_3(\mathrm{c})=\begin{pmatrix}
 			\alpha_1&0&0&\alpha_4 \\
 			0& \beta_2& 1-\alpha_1& 0
 			\end{pmatrix}\simeq\begin{pmatrix}
 			\alpha_1&0&0&a^2\alpha_4 \\
 			0& \beta_2& 1-\alpha_1& 0
 			\end{pmatrix},$  where $ \mathrm{c}=(\alpha_1,\alpha_4,\beta_2)\in\mathbb{F}^3,$ $a\in\mathbb{F}$ and $ a\neq 0$.

  Then we get
$\left\{\begin{array}{rrrr}(\alpha_1-\beta_2)\alpha_4&=&0&\\
\alpha_4(\alpha_1+\beta_2-1)&=&0&\\
 \alpha_4(2\alpha_1-1)&=&0&\\
\beta_2(\alpha_1-\beta_2)&=&0&\\
2\alpha_1^2-3\alpha_1+1&=&0 & \Longleftrightarrow \ \alpha_1=1$ or $\alpha_1=\frac{1}{2}\\
\alpha_4(\alpha_1+\beta_2-1)&=&0&\\
\end{array}\right.$

\textbf{Case 1} $\alpha_1=1:$

$$\left\{\begin{array}{rrrr} \alpha_4(\beta_2-1)&=&0\\
\alpha_4\beta_2&=&0\\
\alpha_4&=&0\\
\beta_2(\beta_2-1)&=&0\\
\alpha_4\beta_2&=&0\\
\end{array}\right. \ \mbox{we get}\  \beta_2(\beta_2-1)=0.$$

\textbf{Case 11} $\beta_2=0$. Then

%a_1 = 1, a_2 = 0, a_3 = 0, a_4 = 0, b_1 = 0, b_2 = 0, b_3 = 0, b_4 = 0

$$\left(\begin{array}{ccccc} 1&0&0&0\\
0&0&0&0
\end{array}\right)$$

\textbf{Case 12} $\beta_2=1$. We get

%a_1 = 1, a_2 = 0, a_3 = 0, a_4 = 0, b_1 = 0, b_2 = 1, b_3 = 0, b_4 = 0
$$\left(\begin{array}{ccccc} 1&0&0&0\\
0&1&0&0
\end{array}\right)$$

\textbf{Case 2} $\alpha_1=\frac{1}{2}:$

One has
$$\left\{\begin{array}{rrr} \alpha_4(2\beta_2-1)&=&0\\
2\beta_2^2-\beta_2&=&0\\
\end{array}\right.$$

%$ \alpha_4*b_2-(1/2)*\alpha_4=0$

\textbf{Case 21} $\alpha_4=0:$ $\Longrightarrow$ \ \ $ 2\beta_2^2-\beta_2=0.$

\textbf{Case 211} $\beta_2=0:$
%a_1 = 1/2, a_2 = 0, a_3 = 0, a_4 = 0, b_1 = 0, b_2 = 0, b_3 = 1/2, b_4 = 0
$$\left(\begin{array}{ccccc} \frac{1}{2}&0&0&0\\
0&0&\frac{1}{2}&0
\end{array}\right)$$

\textbf{Case 212} $\beta_2=\frac{1}{2}:$
%a_1 = 1/2, a_2 = 0, a_3 = 0, b_1 = 0, b_2 = 1/2, b_3 = 1/2, b_4 = 0
$$\left(\begin{array}{ccccc} \frac{1}{2}&0&0&0\\
0&\frac{1}{2}&\frac{1}{2}&0
\end{array}\right)$$

\textbf{Case 22} $a_4\neq 0:$ $\Longrightarrow$ $\beta_2=\frac{1}{2}$

%$\beta_2=1$

%$ -a_4*b_2+(1/2)*a_4= 0$

%$ 2\beta_2^2-\beta_2=0$ %a_1 = 1/2, a_2 = 0, a_3 = 0, b_1 = 0, b_2 = 1/2, b_3 = 1/2, b_4 = 0

$$\left(\begin{array}{ccccc} \frac{1}{2}&0&0&\alpha_4\\
0&\frac{1}{2}&\frac{1}{2}&0
\end{array}\right)\simeq\begin{pmatrix}
 			\frac{1}{2}&0&0&a^2\alpha_4 \\
 			0& \frac{1}{2}& \frac{1}{2}& 0
 			\end{pmatrix},$$  where $ \alpha_4 \in\mathbb{F},$ $a\in\mathbb{F}$ and $ a\neq 0$.

Note that if $\alpha_4\neq 0$ and $\mathbb{F}$ is perfect (particularly, algebraically closed) then $\alpha_4=1$.

For algebras $A_{1}$ and $A_2$ the system of equations (\ref{GSEAs}) also is inconsistent.

Thus we have the following result.

\begin{theor} \label{char0Asso} Any non-trivial $2$-dimensional associative algebra over a field $\mathbb{F},$ $(Char(\mathbb{F})\neq 2,3)$ is isomorphic to only one of the following listed by their matrices of structure constants, such algebras:

\begin{enumerate}
\item  $As_{13}^1:=\left(\begin{array}{ccccc} 0&0&0&0\\
1&0&0&0
\end{array}\right)$
%\end{itemize}
%\item from $A_3$
%\begin{itemize}A3A_a11_a40_b20 := subs([a_1 = 1, a_2 = 0, a_3 = 0, a_4 = 0, b_1 = 0, b_2 = 0, b_3 = 0, b_4 = 0], AX)
\item $As_3^2:=$ $\left(\begin{array}{ccccc} 1&0&0&0\\
0&0&0&0
\end{array}\right)$
\item $As_3^3:=$%A3A_a11_a40_b21 := subs([a_1 = 1, a_2 = 0, a_3 = 0, a_4 = 0, b_1 = 0, b_2 = 1, b_3 = 0, b_4 = 0], AX)
    $\left(\begin{array}{ccccc} 1&0&0&0\\
0&1&0&0
\end{array}\right)$%A3A_a1h_a40_b20 := subs([a_1 = 1/2, a_2 = 0, a_3 = 0, a_4 = 0, b_1 = 0, b_2 = 0, b_3 = 1/2, b_4 = 0], AX)
\item $As_3^4:=$ $\left(\begin{array}{ccccc} \frac{1}{2}&0&0&0\\
0&0&\frac{1}{2}&0
\end{array}\right)$
%A3A_a1h_b2h := subs([a_1 = 1/2, a_2 = 0, a_3 = 0, b_1 = 0, b_2 = 1/2, b_3 = 1/2, b_4 = 0], AX)
\item $As_3^5(\alpha_4):=$ $\left(\begin{array}{ccccc} \frac{1}{2}&0&0&\alpha_4\\
0&\frac{1}{2}&\frac{1}{2}&0
\end{array}\right)\simeq\begin{pmatrix}
 			\frac{1}{2}&0&0&a^2\alpha_4 \\
 			0& \frac{1}{2}& \frac{1}{2}& 0
 			\end{pmatrix},$ where $ \alpha_4 \in\mathbb{F},$ $a\in\mathbb{F}$ and $ a\neq 0$.

%\end{itemize}
\end{enumerate}

\end{theor}

\subsection{Characteristic two}

We apply Theorem \ref{Char2} to verify the algebras given there to be associative.

 All the equations of (\ref{GSEAs}) for algebras $$A_{11,2}(\mathrm{c})=\begin{pmatrix}
 		0&1&1&0 \\
 		\beta_1&0&0&1\\
 		\end{pmatrix}\simeq\begin{pmatrix}
 		0&1&1&0 \\
 		b^2(\beta_1+a^2)&0&0&1\\
 		\end{pmatrix}
 		,\ \mbox{where}\ a, b\in\mathbb{F}, \ b\neq 0$$
 		
 		and $$A_{12,2}=\begin{pmatrix}
 		0&0&0&0 \\
 		1&0&0&0\\
 		\end{pmatrix}
 		$$ in Theorem \ref{Char2} become identities. Therefore, $A_{11,2}$ and $A_{12,2}$ are associative algebras.
%
%$A_{11}$ is associative
%
%
%
%
%
%$A_6$
%
%$a_4=0$
%
%$ a_1= 1$
%%a_2 = 0, a_3 = 0, b_1 = 0, b_2 = 1+a_1, b_3 = a_1, b_4 = 0
%%a_1 = -1, a_2 = 0, a_3 = 0, a_4 = 0, b_1 = 0, b_2 = 0, b_3 = -1, b_4 = 0
%$$\left(\begin{array}{ccccc} 1&0&0&0\\
%0&0&1&0
%\end{array}\right)$$
%
%$A_4$ %a_2 = 1, a_3 = 1, a_4 = 0, b_3 = 1+a_1, b_4 = 1
%
%
%\textbf{Case 1:} $b_2=0$ %[a_1 = 1, a_2 = 0, a_3 = 0, a_4 = 0, b_1 = 0, b_2 = 0, b_3 = 0, b_4 = 0]
%
%$$\left(\begin{array}{ccccc} 1&0&0&0\\
%0&0&0&0
%\end{array}\right)$$
%
%
%\textbf{Case 2:} $b_2\neq 0$ $\Longleftrightarrow$ $a_1=1$
%
%%[a_1 = 1, a_2 = 0, a_3 = 0, a_4 = 0, b_1 = 0, b_2 = 1, b_3 = 0, b_4 = 0]
%$$\left(\begin{array}{ccccc} 1&0&0&0\\
%0&1&0&1
%\end{array}\right)$$
%
%$A_3$ %[a_2 = 0, a_3 = 0, b_1 = 0, b_3 = 1+a_1, b_4 = 0]
%
%$a_4=0$
%
%$ b_2*(1-b_2)=0$ $\Longleftrightarrow$ $b_2=0$ or $b_2=1$
%
%$ a_1=1$
%$$\left(\begin{array}{ccccc} 1&0&0&0\\
%b_1&0&0&1
%\end{array}\right)\ \mbox{where}\ b_1\in \mathbb{F}.$$
%or
%$$\left(\begin{array}{ccccc} 1&0&0&0\\
%b_1&1&0&1
%\end{array}\right)\ \mbox{where}\ b_1\in \mathbb{F}.$$
%

The algebra 
$$A_{4,2}:=\left(\begin{array}{ccccc} 1&1&1&0\\
\beta_1&0&0&1
\end{array}\right)\ \mbox{is associative}.$$

It is easy to see that the algebra $$A_{6,2}:=\left(\begin{array}{ccccc} 1&0&0&0\\
0&0&1&0
\end{array}\right)\ \mbox{also is associative}.$$

Therefore, the following result holds true.

\begin{theor} \label{char2Asso} Any non-trivial $2$-dimensional associative algebra over a field $\mathbb{F},$ $(Char(\mathbb{F})=2)$ is isomorphic to only one of the following listed by their matrices of structure constants, such algebras:

\begin{enumerate}
\item $As_{12,2}^1:=
%\left(\begin{array}{ccccc} 0&0&0&0\\
%\beta_1&0&0&0
%\end{array}\right).$%is associative %a_1 = 0, a_2 = 0, a_3 = 0, a_4 = 0, b_1 = 1, b_2 = 0, b_3 = 0, b_4 = 0
%\begin{itemize}
%\item
\left(\begin{array}{ccccc} 0&0&0&0\\
1&0&0&0
\end{array}\right)$
%\end{itemize}
%a_1 = 0, a_2 = 1, a_3 = 1, a_4 = 0, b_2 = 0, b_3 = 0, b_4 = 1
%\item $As_{11,2}$ is associative %a_1 = 0, a_2 = 0, a_3 = 0, a_4 = 0, b_1 = 1, b_2 = 0, b_3 = 0, b_4 = 0
%\begin{itemize}
\item
    $As_{11,2}^2(\beta_1):=\left(\begin{array}{ccccc} 0&1&1&0\\
\beta_1&0&0&1
\end{array}\right)\ \cong \ \left(\begin{array}{ccccc} 0&1&1&0\\
b^2(\beta_1+a^2)&0&0&1
\end{array}\right), \ \mbox{where}\ a,b \in \mathbb{F} \ \mbox{and} \ b\neq 0.$
%\end{itemize}

%\item $A_{6,2}$ is associative if $\alpha_1=1$ and $a_4=0:$%a_1 = 0, a_2 = 0, a_3 = 0, a_4 = 0, b_1 = 1, b_2 = 0, b_3 = 0, b_4 = 0
%\begin{itemize}
\item
    $As_{6,2}^3:=\left(\begin{array}{ccccc} 1&0&0&0\\
0&0&1&0
\end{array}\right)$
%\end{itemize}

%\item from $A_4$ is associative if $\alpha_1=1\ \mbox{and}\ \beta_2=0:$
%\begin{itemize}%A4_b20_a1n1 := subs([a_1 = -1, a_2 = 1, a_3 = 1, a_4 = 0, b_2 = 0, b_3 = 0, b_4 = 1], AX)
\item  $As_{4,2}^4(\beta_1):=\left(\begin{array}{ccccc} 1&1&1&0\\
\beta_1&0&0&1
\end{array}\right)\ \cong \ \left(\begin{array}{ccccc} 1&1&1&0\\\beta_1+a+a^2&0&0&1
\end{array}\right),\ \mbox{where}\ a, \beta_1 \in \mathbb{F}.$
%b^2(\beta_1+a^2)&0&0&1$
%\item $A_6^3:=$ $\left(\begin{array}{ccccc} \frac{1}{2}&0&0&0\\
%0&0&\frac{1}{2}&0
%\end{array}\right)$
%\end{itemize}
%\item from $A_3$ the associative cases are:
%\begin{itemize}%A3_a11_a40_b20 := subs([a_1 = 1, a_2 = 0, a_3 = 0, a_4 = 0, b_1 = 0, b_2 = 0, b_3 = 0, b_4 = 0], AX)
%\item  $\alpha_1=1, \alpha_4=0\ \mbox{and}\ \beta_2=0:$
\item $As_{3,2}^5:=$ $\left(\begin{array}{ccccc} 1&0&0&0\\
0&0&0&0
\end{array}\right),$
%A3_a11_a40_b21 := subs([a_1 = 1, a_2 = 0, a_3 = 0, a_4 = 0, b_1 = 0, b_2 = 1, b_3 = 0, b_4 = 0], AX)
\item $As_{3,2}^6:=$ $\left(\begin{array}{ccccc} 1&0&0&0\\
0&1&0&0
\end{array}\right)$
%\end{itemize}
\end{enumerate}

\end{theor}

%a_1 = 0, a_2 = 1, a_3 = 1, a_4 = 0, b_1 = 1, b_2 = 0, b_3 = 1, b_4 = -1

%LJM dagi maqolada $(Char(\mathbb{F})=2)$ assotiativ deb berilgan quyidagi algebra menda assotiativ chiqmayapti:
%
%$$A_{3,2}(\beta_1=1,\beta_2=0):=\left(\begin{array}{ccccc} 0&1&1&0\\
%1&0&1&-1
%\end{array}\right)$$

\subsection{Characteristic three} \label{char3Asso}

In this case the associative algebras come out from the following classes of Theorem \ref{Char3}.

It is immediate to get that the algebra $A_{13,3}$ is associative. In these case all the equations of the system (\ref{GSEAs}) turn into identities.

%$$A_{13,3}:=\left(\begin{array}{ccccc} 0&0&0&0\\
%1&0&0&0
%\end{array}\right)\ \mbox{is associative}$$

%$A_{10}$ %a_1 = 0, a_2 = 0, a_3 = 0, a_4 = 1, b_2 = 0, b_3 = 0, b_4 = 0

%$$As_{10,3}:=\left(\begin{array}{ccccc} 0&0&0&1\\
%0&0&0&0
%\end{array}\right)\ \mbox{also is associative}$$

Let us consider
$$A_{3,3}(\mathrm{c})=\begin{pmatrix}
 		\alpha_1&0&0&\alpha_4 \\
 		0& \beta_2& 1-\alpha_1& 0
 		\end{pmatrix}\cong \begin{pmatrix}
 		\alpha_1&0&0&a^2\alpha_4 \\
 		0& \beta_2& 1-\alpha_1& 0
 		\end{pmatrix},$$ where $ \mathrm{c}=(\alpha_1,\alpha_4,\beta_2)\in\mathbb{F}^3,$ $ a\in\mathbb{F}$ and $ a\neq 0$.

 The system of equations (\ref{GSEAs}) is equivalent to

%$(\alpha_1-\beta_2)\alpha_4=0$
%
%$\alpha_4(\alpha_1+\beta_2-1)=0$
%
%$ \alpha_4(2\alpha_1-1)=0$
%
%$\beta_2(\alpha_1-\beta_2)=0$
%
%$2\alpha_1^2-3\alpha_1+1=0$
%
%$\alpha_4(\alpha_1+\beta_2-1)$

\begin{equation} \label{A3}
\left\{\begin{array}{lll}
 (\alpha_1-\beta_2)\alpha_4&=&0\\   \alpha_4(\alpha_1+\beta_2-1)&=&0\\ \alpha_4(2\alpha_1-1)&=&0\\ \beta_2(\alpha_1-\beta_2)&=&0\\  2\alpha_1^2-3\alpha_1+1&=&0\\  \alpha_4(\alpha_1+\beta_2-1)&=&0\\
\end{array}\right.
\end{equation}

From (\ref{A3}) one has $2\alpha_1^2-3\alpha_1+1=0$ $\Longleftrightarrow$ $\alpha_1=1$ or $\alpha_1=2.$

\textbf{Case 1:} $\alpha_1=1$ Then (\ref{A3}) is equivalent to $ \beta_2^2-\beta_2=0$. Therefore, we have two subcases:

%$a_4*(b_2-1)=0$
%
%$a_4*b_2=0$
%
%$ a_4=0$
%
%$ b_2^2-b_2=0$
%\begin{equation}
%\begin{array}{lll}
%a_4(b_2-1)&=&0\\ a_4b_2&=&0\\ a_4&=&0\\ b_2^2-b_2&=&0\\
%\end{array}
%\end{equation}

\textbf{Case 11.} Let $\beta_2=0.$ Then we get

%a_1 = 1, a_2 = 0, a_3 = 0, a_4 = 0, b_1 = 0, b_2 = 0, b_3 = 0, b_4 = 0
$$A_{3,3}:=\left(\begin{array}{ccccc} 1&0&0&0\\
0&0&0&0
\end{array}\right)\ \mbox{is associative}.$$

\textbf{Case 12:} Let $\beta_2=1$. Then one obtains that %a_1 = 1, a_2 = 0, a_3 = 0, a_4 = 0, b_1 = 0, b_2 = 1, b_3 = 0, b_4 = 0
$$A_{3,3}:=\left(\begin{array}{ccccc} 1&0&0&0\\
0&1&0&0
\end{array}\right)\ \mbox{is associative}.$$

\textbf{Case 2:} If $\alpha_1=2$ then (\ref{A3}) is equivalent to $ \beta_2^2-\beta_2=0$. Considering two subcases for $\beta_2=0$ (which implies $\alpha_4=0$) and $\beta_2=2$ we obtain the following two associative algebras:

%$(\alpha_1-\beta_2)=0$
%
%$(\alpha_1+\beta_2-1)=0$
%
%$ (2\alpha_1-1)=0$ $\Longleftrightarrow$ $\alpha_1=2$
%
%$\beta_2(\alpha_1-\beta_2)=0$
%
%$2\alpha_1^2+1=0$%a_1 = 2, a_2 = 0, a_3 = 0, b_1 = 0, b_2 = 2, b_3 = 2, b_4 = 0
%
%\begin{equation}
%\begin{array}{lll}
%a_4(b_2-2)&=&0\\ -a_4b_2+2a_4&=&0\\  -b_2^2+2b_2&=&0\\
%\end{array}
%\end{equation}

%\textbf{Case 21.} Let $b_2=0.$ This implies $a_4=0$. Therefore, we get
%
$$A_{3,3}:=\left(\begin{array}{ccccc} 2&0&0&0\\
0&0&2&0
\end{array}\right)$$
and
%\textbf{Case 22.} Let now $b_2=2.$ Then, one gets
%a_1 = 2, a_2 = 0, a_3 = 0, b_1 = 0, b_2 = 2, b_3 = 2, b_4 = 0
$$A_{3,3}:=\left(\begin{array}{ccccc} 2&0&0&\alpha_4\\
0&2&2&0
\end{array}\right)\cong \begin{pmatrix}
 		2&0&0&a^2\alpha_4 \\
 		0& 2& 2& 0
 		\end{pmatrix},$$ where $ \alpha_4 \in\mathbb{F},$ $ a\in\mathbb{F}$ and $ a\neq 0$.

There are no associative algebras generated from the other classes of Theorem \ref{Char3}.

Thus, we have the following theorem.

\begin{theor} \label{char3Asso} Any non-trivial $2$-dimensional associative algebra over a field $\mathbb{F},$ $(Char(\mathbb{F})=3)$ is isomorphic to only one of the following listed by their matrices of structure constants, such algebras:

\begin{enumerate}
\item $As_{13,3}^1:=\left(\begin{array}{ccccc} 0&0&0&0\\
1&0&0&0
\end{array}\right)$.
%\itemA3_a11_a40_b20 := subs([a_1 = 1, a_2 = 0, a_3 = 0, a_4 = 0, b_1 = 0, b_2 = 0, b_3 = 0, b_4 = 0], AX)
\item
    $As_{3,3}^2:=\left(\begin{array}{ccccc} 1&0&0&0\\
0&0&0&0
\end{array}\right)$
%\end{itemize}
%A3_a11a40b21 := subs([a_1 = 1, a_2 = 0, a_3 = 0, a_4 = 0, b_1 = 0, b_2 = 1, b_3 = 0, b_4 = 0], AX)
\item
    $As_{3,3}^3:=\left(\begin{array}{ccccc} 1&0&0&0\\
0&1&0&0
\end{array}\right)$
\item%A3_a12_a40_b20 := subs([a_1 = 2, a_2 = 0, a_3 = 0, a_4 = 0, b_1 = 0, b_2 = 0, b_3 = 2, b_4 = 0], AX)
    $As_{3,3}^4:=\left(\begin{array}{ccccc} 2&0&0&0\\
0&0&2&0
\end{array}\right)$
%A3_a12_a40_b22 := subs([a_1 = 2, a_2 = 0, a_3 = 0, b_1 = 0, b_2 = 2, b_3 = 2, b_4 = 0], AX)
\item $As_{3,3}^5(\alpha_4):=\left(\begin{array}{ccccc} 2&0&0&\alpha_4\\
0&2&2&0
\end{array}\right)\cong \begin{pmatrix}
 		2&0&0&a^2\alpha_4 \\
 		0& 2& 2& 0
 		\end{pmatrix},$ where $ \alpha_4 \in\mathbb{F},$ $ a\in\mathbb{F}$ and $ a\neq 0$.
\end{enumerate}

\end{theor}

\section{Automorphism groups}

In this section we describe the automorphism groups of algebras from Theorems \ref{char0Asso}, \ref{char2Asso} and \ref{char3Asso}. The author believes such automorphism groups can be obtained easily. But the lists of associative algebras in the theorems are over arbitrary field and we do it here for the paper to be self-contained. We need the automorphism groups in the next section to verify whether some of two-dimensional diassociative algebras found there isomorphic or not. Let $g=\left(\begin{array}{lll} x & y\\ z & t \end{array}\right)$ with $xt\neq yz$. The equation (\ref{1}) is equivalent to

\begin{equation} \label{SEAut3}
\left\{\begin{array}{lll}
\alpha_1x^2+((\alpha_2+\alpha_3)z+\alpha_1)x+\alpha_4z^2-\beta_1y&=&0\\
(\alpha_1y+\alpha_2(t-1))x+(\alpha_3z-\beta_2)y+\alpha_4tz&=&0\\
 (\alpha_1y+\alpha_3(t-1))x+(\alpha_2z-\beta_3)y+\alpha_4tz&=&0\\
  \alpha_1y^2+((\alpha_2+\alpha_3)t-\beta_4)y+\alpha_4(t^2-x)&=&0\\
  \beta_4z^2+((\beta_2+\beta_3)x-\alpha_1)z+\beta_1(x^2-t)&=&0\\
  (\beta_4z+\beta_2(x-1))t+(\beta_3y-\alpha_2)z+\beta_1xy&=&0\\
   (\beta_4z+\beta_3(x-1))t+(\beta_2y-\alpha_3)z+\beta_1xy&=&0\\
   \beta_4t^2+((\beta_2+\beta_3)y-\beta_4)t+\beta_1y^2-\alpha_4z&=&0\\
\end{array}\right.
\end{equation}

\subsection{Characteristic of $\mathbb{F}$ is not 2 and 3}\emph{}\\

For $As_{13}^1$ the system (\ref{SEAut3}) is equivalent to
$\left\{\begin{array}{rrr} y&=&0\\ x^2-t&=&0\\
\end{array}\right.$

 Therefore,
$$Aut(As_{13}^1)=Aut\left(\begin{pmatrix}
 			0&0&0&0 \\
 			1&0&0&0
 			\end{pmatrix}\right)=\left\{\left(\begin{array}{lll} x & 0\\ z & x^2 \end{array}\right)\big| \ x\neq 0\right\}.$$

Consider $As_3^2$. Then as the system (\ref{SEAut3}) we get

$$\left\{\begin{array}{rrr} x(x-1)&=&0\\
xy&=&0\\ y&=&0\\ z&=&0\\ \end{array}\right.$$

Hence,

$$Aut(As_{3}^2)=Aut\left(\begin{pmatrix}
 			1&0&0&0 \\
 			0&0&0&0
 			\end{pmatrix}\right)=\left\{\left(\begin{array}{lll} 1 & 0\\ 0 & t \end{array}\right)\big| \ t\neq 0\right\}.$$

Consider $As_3^3$. Then
$$\left\{\begin{array}{rrr}x(x-1)&=&0\\ y&=&0\\ z(x-1)&=&0\\ t(x-1)&=&0\\ \end{array}\right.$$

and

$$Aut(As_{3}^3)=Aut\left(\begin{pmatrix}
 			1&0&0&0 \\
 			0&1&0&0
 			\end{pmatrix}\right)=\left\{\left(\begin{array}{lll} 1 & 0\\ z & t \end{array}\right)\big| \ t\neq 0\right\}.$$

Consider $As_3^4:=\begin{pmatrix}
 			\frac{1}{2}&0&0&0 \\
 			0&0&\frac{1}{2}&0
 			\end{pmatrix}$. We get
$$\left\{\begin{array}{rrr} x-x^2&=&0\\ y&=&0\\ z(x-1)&=&0\\  t(x-1)&=&0 \\ \end{array}\right.$$

  and

 $$Aut(As_{3}^4)=Aut\left(\begin{pmatrix}
 			\frac{1}{2}&0&0&0 \\
 			0&0&\frac{1}{2}&0
 			\end{pmatrix}\right)=\left\{\left(\begin{array}{lll} 1 & 0\\ z & t \end{array}\right)\big| \ t\neq 0\right\}.$$

Let us now consider $As_3^5:=\begin{pmatrix}
 			\frac{1}{2}&0&0&\alpha_4 \\
 			0&\frac{1}{2}&\frac{1}{2}&0
 			\end{pmatrix}$. Then
$$\left\{\begin{array}{rrr} x-x^2-2\alpha_4z^2&=&0\\ (x-1)y+2\alpha_4zt&=&0\\  2x\alpha_4-y^2-2\alpha_4t^2&=& 0 \\ z-2xz&=&0\\ (x-1)t+zy&=&0\\ \alpha_4z-ty &=&0\end{array}\right.$$

 %$\left(\begin{array}{lrr} 1 & 0\\ 0 & \pm 1 \end{array}\right)$
 The solution to the system is $\left\{\begin{array}{lll}
 \{x=1, y=0, z=0, t \ \mbox{is any non-zero}\}&\ \mbox{if}\ \alpha_4=0&\\  \{x=1, y=0, z=0, t=\pm 1\}& \ \mbox{if}\ \alpha_4\neq 0&\\
\end{array}\right.$

i.e.,
$$Aut(As_{3}^5(0))=Aut\left(\begin{pmatrix}
 			\frac{1}{2}&0&0&0 \\
 			0&\frac{1}{2}&\frac{1}{2}&0
 			\end{pmatrix}\right)=\left\{\left(\begin{array}{lrr} 1 & 0\\ 0 & t \end{array}\right)\big| \ t\neq 0 \right\},$$
$$Aut(As_{3}^{5}(\alpha_4))=Aut\left(\begin{pmatrix}
 			\frac{1}{2}&0&0&0 \\
 			0&\frac{1}{2}&\frac{1}{2}&0
 			\end{pmatrix}\right)=\left\{ I=\left(\begin{array}{crr} 1 & 0\\ 0 & \pm1 \end{array}\right)\right\}.$$

\subsection{Characteristic of $\mathbb{F}$ is two}\emph{}\\

Consider $As_{12,2}^1=\begin{pmatrix}
 		0&0&0&0 \\
 		1&0&0&0\\
 		\end{pmatrix}
 		$. From (\ref{SEAut3}) we get $y=0$ and $t=x^2$. Therefore,
$$Aut(As_{12,2}^1)=\left\{ \left(\begin{array}{lrr} x & 0\\ z & x^2 \end{array}\right)\big| \ \mbox{where}\ x\neq 0, z \in \mathbb{F}\right\}.$$

Let us take $As_{11,2}^2(\beta_1)=\begin{pmatrix}
 		0&1&1&0 \\
 		\beta_1&0&0&1\\
 		\end{pmatrix}\simeq\begin{pmatrix}
 		0&1&1&0 \\
 		b^2(\beta_1+a^2)&0&0&1\\
 		\end{pmatrix}
 		,$ where $a, b\in\mathbb{F}$ and $ b\neq 0$. Then we get

$$\left\{\begin{array}{rrr} t&=&1\\ y&=&0\\ z&=&\beta_1(x-1)\\ \end{array} \right.$$ and
$$Aut(As_{11,2}^2)=\left\{ \left(\begin{array}{crr} x & 0\\ \beta_1(x-1) & 1 \end{array}\right)\big| \ \mbox{where}\ x\neq 0 \in \mathbb{F}\right\}.$$

Consider $As_{6,2}^3=\begin{pmatrix}%A6_a1n_a40 := subs([a_1 = -1, a_2 = 0, a_3 = 0, a_4 = 0, b_1 = 0, b_2 = 0, b_3 = -1, b_4 = 0], AX)
 		1&0&0&0 \\
 		0& 0& 1& 0
 		\end{pmatrix}$.

 Then (\ref{SEAut3}) is equivalent to
$$\left\{\begin{array}{rrr} x^2-x&=&0\\ xy&=&0\\ y(x-1)&=&0\\ y^2+y&=&0\\ z(x-z-1)&=&0\\ z(t-y)&=&0\\ t(x-z-1)&=&0\\ t(t-y-1)&=&0\\
\end{array}\right.$$ and
$$Aut(As_{6,2}^3)=\left\{ I=\left(\begin{array}{crr} 1 & 0\\ 0 & 1 \end{array}\right)\right\}.$$

Consider $ As_{4,2}^4(\beta_1):=\begin{pmatrix}
 		1&1&1&0 \\
 		\beta_1& 0& 0& 1
 		\end{pmatrix}\simeq\begin{pmatrix}
 		1&1&1&0 \\
 		\beta_1+a+a^2& 0& 0& 1
 		\end{pmatrix}
 		$.
%A4_b20_a1n1 := subs([a_1 = -1, a_2 = 1, a_3 = 1, a_4 = 0, b_2 = 0, b_3 = 0, b_4 = 1], AX)
 Then from the system of equations (\ref{SEAut3}) we obtain
$$\left\{\begin{array}{rrr} \beta_1y&=&0\\ (t+y-1)x-zy&=&0\\ y(y-1)&=&0\\ (x-t)b_1&=&0\\ tz+z&=&0.\\
\end{array}\right. $$

\begin{itemize}
\item if $\beta_1=0$ we get
    $$Aut(As_{4,2}^4(0))=\left\{ \left(\begin{array}{crr} x & 0\\ z & 1 \end{array}\right)\big| \ x \neq 0, \ z \in \mathbb{F} \right\}.$$
\item if $\beta_1\neq 0$ then
$$Aut(As_{4,2}^4(\beta_1))=\left\{ \left(\begin{array}{crr} 1 & 0\\ z & 1 \end{array}\right)\big| \ z \in \mathbb{F} \right\}.$$
\end{itemize}
%A3_a11_a40_b20 := subs([a_1 = 1, a_2 = 0, a_3 = 0, a_4 = 0, b_1 = 0, b_2 = 0, b_3 = 0, b_4 = 0], AX)
Consider $ As_{3,2}^5:=\begin{pmatrix}
 		1&0&0&0 \\
 		0& 0& 0& 0
 		\end{pmatrix}$.

 Then (\ref{SEAut3}) becomes
$$\left\{\begin{array}{rrr} x^2+x&=&0\\ xy&=&0\\ y&=&0\\ z&=&0.\\ \end{array}\right. $$

Hence,
$$Aut(As_{3,2}^5)=\left\{ \left(\begin{array}{crr} 1 & 0\\ 0 & t \end{array}\right)\big| \ t \neq 0 \in \mathbb{F}\right\}.$$
Consider
%A3_a11_a40_b21 := subs([a_1 = 1, a_2 = 0, a_3 = 0, a_4 = 0, b_1 = 0, b_2 = 1, b_3 = 0, b_4 = 0], AX)
$ As_{3,2}^6:=\begin{pmatrix}
 		1&0&0&0 \\
 		0& 1& 0& 0
 		\end{pmatrix}$. Then
$$\left\{\begin{array}{rrr} y&=&0\\ xz+z&=&0\\ tx+t&=&0\end{array}\right. $$ Therefore,
$$Aut(As_{3,2}^6)=\left\{ \left(\begin{array}{crr} 1 & 0\\ z & t \end{array}\right)\big| \ z \in \mathbb{F} \ \mbox{and}\ t \neq 0 \in \mathbb{F}\right\}.$$

\subsection{Characteristic of $\mathbb{F}$ is three}\emph{}

Consider $As_{13,3}^1:=\begin{pmatrix}
 		1&0&0&0 \\
 		1& 0& 0& 0
 		\end{pmatrix}$.

$$\left\{\begin{array}{rrr} y&=&0\\ x^2+t&=&0 \end{array}\right. $$

Therefore,

$$Aut(As_{13,3}^1)=\left\{ \left(\begin{array}{crr} x & 0\\ z & 2x^2 \end{array}\right)\right\}.$$

Consider $As_{3,3}^2:=\begin{pmatrix}
 		1&0&0&0 \\
 		0& 0& 0& 0
 		\end{pmatrix}$.
 
 From (\ref{SEAut3}) we obtain
$$\left\{\begin{array}{rrr} x(x-1)&=&0\\ y&=&0\\ z&=&0\end{array}\right. $$

Hence,
$$Aut(As_{3,3}^2)=\left\{ \left(\begin{array}{crr} 1 & 0\\ z & t \end{array}\right)\big| \ t\neq 0\right\}.$$

Consider $As_{3,3}^3:=\begin{pmatrix}
 		1&0&0&0 \\
 		0& 1& 0& 0
 		\end{pmatrix}$.
Then
%A3_a11a40b21 := subs([a_1 = 1, a_2 = 0, a_3 = 0, a_4 = 0, b_1 = 0, b_2 = 1, b_3 = 0, b_4 = 0], AX)
$$\left\{\begin{array}{rrr} x(x-1)&=&0\\ y&=&0\\ (x+1)z&=&0\\ t(x-1)&=&0\end{array}\right. $$
and
$$Aut(As_{3,3}^3)=\left\{ \left(\begin{array}{crr} 1 & 0\\ 0 & t \end{array}\right)\big| \ t\neq 0\right\}.$$

%A3_a12_a40_b20 := subs([a_1 = 2, a_2 = 0, a_3 = 0, a_4 = 0, b_1 = 0, b_2 = 0, b_3 = 2, b_4 = 0], AX)
If $A=As_{3,3}^4:=\begin{pmatrix}
 		2&0&0&0 \\
 		0& 0& 2& 0
 		\end{pmatrix}$ then (\ref{SEAut3}) implies
$$\left\{\begin{array}{rrr} x^2-x&=&0\\ y&=&0\\ (2x-1)z-x^2+t&=&0\\ (x-1)t&=&0\end{array}\right. $$

Therefore,

$$Aut(As_{3,3}^4)=\left\{ \left(\begin{array}{crr} 1 & 0\\ 1+2t & t \end{array}\right)\big| \ t \neq 0 \right\}.$$

The system of equations (\ref{SEAut3}) for the group of automorphisms of $As_{3,3}^5(\alpha_4):=\begin{pmatrix}
 		2&0&0&\alpha_4 \\
 		0& 2& 2& 0
 		\end{pmatrix}$ is
%A3_a12_a40_b22 := subs([a_1 = 2, a_2 = 0, a_3 = 0, b_1 = 0, b_2 = 2, b_3 = 2, b_4 = 0], AX)
%$ (-2x+2)y-\alpha_4zt=0$
%
%$ (-t^2+x)\alpha_4-2y^2=0$
%
%$2z(x+1)=0$
%
%$ (-2x+2)t-2zy=0$
%
%$ \alpha_4z+2ty=0$
%
%Simplify:

\begin{equation}
\left\{\begin{array}{lll}
\alpha_4z^2+2x^2+x+y&=&0\\ \alpha_4tz+2xy&=&0\\ (2x+1)y+\alpha_4zt&=&0\\ (t^2-x)\alpha_4+2y^2&=&0\\ (2x+1)z+2x^2+t&=&0\\ y(x+z)&=&0\\ (2x+1)t+2xy&=&0\\ \alpha_4z+ty+y^2&=&0
\end{array}\right.
\end{equation}

The solution to the system is $\left\{\begin{array}{lll}
 \{x=1, y=0, z=0, t=1\}&\ \mbox{if}\ \alpha_4=0&\\  \{x=1, y=0, z\ \mbox{is any,}\ t=1\}& \ \mbox{if}\ \alpha_4\neq 0&\\
\end{array}\right.$

Thus,
$$Aut(As_{3,3}^{4,1}(0))=Aut\left(\left(\begin{array}{ccccc} 2&0&0&0\\
0&2&2&0
\end{array}\right)\right)=\left\{ I=\left(\begin{array}{crr} 1 & 0\\ 0 & 1 \end{array}\right)\right\}$$

$$Aut(As_{3,3}^{4,2}(\alpha_4))=Aut\left(\left(\begin{array}{ccccc} 2&0&0&\alpha_4\\
0&2&2&0
\end{array}\right)\right)=\left\{ \left(\begin{array}{crr} 1 & 0\\ z & 1 \end{array}\right),\ z \in \mathbb{F} \right\}.$$

% Consider $$As_{3,3}^5(\alpha_4):=\left(\begin{array}{ccccc} 2&0&0&\alpha_4\\
%0&2&2&0
%\end{array}\right)$$
%\begin{equation} \label{GSE}
%\begin{array}{lll}
%-a_4z^2-2x^2+2x&=&0\\ (-2x+2)y-a_4zt&=&0\\ (-2x+2)y-a_4zt&=&0\\ (-t^2+x)a_4-2y^2&=&0\\ -4xz+2z&=&0\\ (-2x+2)t-2zy&=&0\\ (-2x+2)t-2zy&=&0\\ a_4z-4ty&=&0\\
%\end{array}
%\end{equation}

\section{Classification of two-dimensional associative dialgebras}

In this section we classify all two-dimensional associative dialgebras over any basic field. As was mentioned earlier
a di-algebra can be given by  two $2\times 4$ matrices \[A=\left(\begin{array}{cccc} \alpha_1 & \alpha_2 & \alpha_3 &\alpha_4\\ \beta_1 & \beta_2 & \beta_3 &\beta_4\end{array}\right)\ \mbox{and}\ B=\left(\begin{array}{ccccc} \gamma_1&\gamma_2&\gamma_3&\gamma_4\\
\delta_1&\delta_2&\delta_3&\delta_4
\end{array}\right)\] corresponding to the binary operations $\dashv$ and $\vdash$, respectively. The matrix equations  (\ref{DAA3}) in terms of entries of $A$ and $B$ can be written as follows

AXIOM 1: $$A(A \otimes I)-A(I \otimes A)=0$$
\begin{equation} \label{GSE}
\begin{array}{lll}
  \beta_1(\alpha_2-\alpha_3)&=&0\\
  \alpha_2\beta_2-\alpha_4\beta_1&=&0\\
 (\alpha_1-\beta_3)\alpha_2-\alpha_3(\alpha_1-\beta_2)&=&0\\
 (\alpha_1-\beta_2)\alpha_4-\alpha_2(\alpha_2-\beta_4)&=&0\\
  \alpha_3\beta_3-\alpha_4\beta_1&=&0\\
  \alpha_4(\beta_2-\beta_3)&=&0\\
 (\alpha_1-\beta_3)\alpha_4-\alpha_3(\alpha_3-\beta_4)&=&0\\
 \alpha_4(\alpha_2-\alpha_3)&=&0\\
 \beta_1(\beta_2-\beta_3)&=&0\\
 (\alpha_2-\beta_4)\beta_1-\beta_2(\alpha_1-\beta_2)&=&0\\
% \beta_1(\alpha_2-\alpha_3)&=&0\\
% \alpha_2\beta_2-\alpha_4\beta_1&=&0\\
 (\alpha_3-\beta_4)\beta_1-\beta_3(\alpha_1-\beta_3)&=&0\\
 (\alpha_3-\beta_4)\beta_2-\beta_3(\alpha_2-\beta_4)&=&0\\
% \alpha_3\beta_3-\alpha_4\beta_1&=&0\\
% \alpha_4(\beta_2-\beta_3)&=&0\\
%\begin{equation} \label{GSE}
%\begin{array}{lll}
\end{array}
\end{equation}

AXIOM 2: $$A(I \otimes A)-A(I \otimes B)=0$$
\begin{equation}
\begin{array}{lll}
\alpha_1^2-\alpha_1\gamma_1+\alpha_2\beta_1-\alpha_2\delta_1&=&0\\ \alpha_1\alpha_2-\alpha_1\gamma_2+\alpha_2\beta_2-\alpha_2\delta_2&=&0\\ \alpha_1\alpha_3-\alpha_1\gamma_3+\alpha_2\beta_3-\alpha_2\delta_3&=&0\\ \alpha_1\alpha_4-\alpha_1\gamma_4+\alpha_2\beta_4-\alpha_2\delta_4&=&0\\ \alpha_1\alpha_3-\alpha_3\gamma_1+\alpha_4\beta_1-\alpha_4\delta_1&=&0\\ \alpha_2\alpha_3-\alpha_3\gamma_2+\alpha_4\beta_2-\alpha_4\delta_2&=&0\\ \alpha_3^2-\alpha_3\gamma_3+\alpha_4\beta_3-\alpha_4\delta_3&=&0\\ \alpha_3\alpha_4-\alpha_3\gamma_4+\alpha_4\beta_4-\alpha_4\delta_4&=&0\\\alpha_1\beta_1+\beta_1\beta_2-\beta_1\gamma_1-\beta_2\delta_1&=&0\\ \alpha_2\beta_1-\beta_1\gamma_2+\beta_2^2-\beta_2\delta_2&=&0\\ \alpha_3\beta_1-\beta_1\gamma_3+\beta_2\beta_3-\beta_2\delta_3&=&0\\ \alpha_4\beta_1-\beta_1\gamma_4+\beta_2\beta_4-\beta_2\delta_4&=&0\\ \alpha_1\beta_3+\beta_1\beta_4-\beta_3\gamma_1-\beta_4\delta_1&=&0\\ \alpha_2\beta_3+\beta_2\beta_4-\beta_3\gamma_2-\beta_4\delta_2&=&0\\ \alpha_3\beta_3+\beta_3\beta_4-\beta_3\gamma_3-\beta_4\delta_3&=&0\\ \alpha_4\beta_3-\beta_3\gamma_4+\beta_4^2-\beta_4\delta_4&=&0\\
\end{array}
\end{equation}

AXIOM 3: $$A(B\otimes I)-B(I \otimes A)=0$$
\begin{equation}
\begin{array}{lll}
\alpha_3\delta_1-\beta_1\gamma_2&=&0\\
 \alpha_4\delta_1-\beta_2\gamma_2&=&0\\
  (\delta_2-\gamma_1)\alpha_3+\gamma_2(\alpha_1-\beta_3)&=&0\\
   (\delta_2-\gamma_1)\alpha_4+\gamma_2(\alpha_2-\beta_4)&=&0\\
    \alpha_3\delta_3-\beta_1\gamma_4&=&0\\
     \alpha_4\delta_3-\beta_2\gamma_4&=&0\\
      (\gamma_3-\delta_4)\alpha_3-\gamma_4(\alpha_1-\beta_3)&=&0\\ (\gamma_3-\delta_4)\alpha_4-\gamma_4(\alpha_2-\beta_4)&=&0\\
      (\gamma_1-\delta_2)\beta_1-\delta_1(\alpha_1-\beta_3)&=&0\\ \beta_2(\gamma_1-\delta_2)-\delta_1(\alpha_2-\beta_4)&=&0\\
      \beta_1(\gamma_3-\delta_4)-\delta_3(\alpha_1-\beta_3)&=&0\\ (\gamma_3-\delta_4)\beta_2-\delta_3(\alpha_2-\beta_4)&=&0\\
       \end{array}
\end{equation}

AXIOM 4: $$B(A \otimes I)-B(B \otimes I)=0$$
\begin{equation}
\begin{array}{lll}
\alpha_1\gamma_1-\gamma_1^2+(-\delta_1+\beta_1)\gamma_3&=&0\\
 (\beta_1-\delta_1)\gamma_4+\gamma_2(\alpha_1-\gamma_1)&=&0\\
  \gamma_1(\alpha_2-\gamma_2)+\gamma_3(\beta_2-\delta_2)&=&0\\
   \alpha_2\gamma_2-\gamma_2^2+\gamma_4(\beta_2-\delta_2)&=&0\\ (\beta_3-\gamma_1-\delta_3)\gamma_3+\alpha_3\gamma_1&=&0\\ \gamma_2(\alpha_3-\gamma_3)+\gamma_4(\beta_3-\delta_3)&=&0\\ (\alpha_4-\gamma_4)\gamma_1+\gamma_3(\beta_4-\delta_4)&=&0\\ (\beta_4-\gamma_2-\delta_4)\gamma_4+\alpha_4\gamma_2&=&0\\(\alpha_1-\delta_3-\gamma_1)\delta_1+\beta_1\delta_3&=&0\\ (\alpha_1-\gamma_1)\delta_2+\delta_4(-\delta_1+\beta_1)&=&0\\ (\alpha_2-\gamma_2)\delta_1+\delta_3(\beta_2-\delta_2)&=&0\\ (\alpha_2-\gamma_2-\delta_4)\delta_2+\beta_2\delta_4&=&0\\ (\alpha_3-\gamma_3)\delta_1+\delta_3(\beta_3-\delta_3)&=&0\\ (\alpha_3-\gamma_3)\delta_2+\delta_4(\beta_3-\delta_3)&=&0\\ (\alpha_4-\gamma_4)\delta_1+\delta_3(\beta_4-\delta_4)&=&0\\ (\alpha_4-\gamma_4)\delta_2+\delta_4(\beta_4-\delta_4)&=&0\\
\end{array}
\end{equation}

AXIOM 5: $$B(B \otimes I)-B  (B \otimes I)=0$$
\begin{equation}
\begin{array}{lll}
\delta_1(\gamma_2-\gamma_3)&=&0\\
\gamma_2\delta_2-\gamma_4\delta_1&=&0\\
 \gamma_1(\gamma_2-\gamma_3)-\gamma_2\delta_3+\gamma_3\delta_2&=&0\\ \gamma_2^2-\gamma_2\delta_4-\gamma_4(\gamma_1-\delta_2)&=&0\\
  \gamma_3\delta_3-\gamma_4\delta_1&=&0\\
   \gamma_4(\delta_2-\delta_3)&=&0\\
   \gamma_3\delta_4-\gamma_3^2+\gamma_4(\gamma_1-\delta_3)&=&0\\ \gamma_4(\gamma_2-\gamma_3)&=&0\\
   \delta_1(\delta_2-\delta_3)&=&0\\
    (\gamma_2-\delta_4)\delta_1-\delta_2(\gamma_1-\delta_2)&=&0\\
    (\gamma_3-\delta_4)\delta_1-\delta_3(\gamma_1-\delta_3)&=&0\\ (\gamma_3-\delta_4)\delta_2-\delta_3(\gamma_2-\delta_4)&=&0\\

\end{array}
\end{equation}

For $A$ we take MSC of Theorems \ref{char0Asso}, \ref{char2Asso}, \ref{char3Asso} for a basic field is not characteristic 2,3, characteristic 2 and characteristic 3 cases, respectively. The entries of $B$ we consider as unknowns: $$\left(\begin{array}{ccccc} \gamma_1&\gamma_2&\gamma_3&\gamma_4\\
\delta_1&\delta_2&\delta_3&\delta_4
\end{array}\right)$$ Substitute these $A$ and $B$ into the matrix equations  (\ref{DAA3}) to get the system of equation MSC chosen $A$ with unknown entries of $B$. Solving the system of equations we get a diassociative algebra generated by $A$. Acting by the automorphism group of $A$ we verify whether the generated by $A$ diassociative algebras are isomorphic or not.

\subsection{Characteristic of $\mathbb{F}$ is not two and three}\emph{}

Consider $As_{13}^1$.

AXIOM 1 and AXIOM 2 $\Longrightarrow$ $\gamma_1=\gamma_2=\gamma_3=\gamma_4=0$

AXIOM 3 $\Longrightarrow$ $\delta_2=\delta_4=0$

AXIOM 4 gives \ \
$\left\{\begin{array}{rrr}\delta_1\delta_3-\delta_3&=&0\\ \delta_3&=&0\end{array}\right.$

AXIOM 5 holds true.

%a_1 = 0, a_2 = 0, a_3 = 0, a_4 = 0, b_1 = 1, b_2 = 0, b_3 = 0, b_4 = 0, c_1 = 0, c_2 = 0, c_3 = 0, c_4 = 0, d_2 = 0, d_3 = 0, d_4 = 0

Therefore,

$$D_{13}^1:=\left\{A=\left(\begin{array}{cccc} 0 & 0 & 0 &0\\ 1 & 0 & 0 &0\end{array}\right),\ \ B=\left(\begin{array}{cccc} 0 & 0 & 0 &0\\ \delta_1& 0& 0 & 0\end{array}\right)\right\}$$

is a diassociative algebra generated by $As_{13}^1$.

Consider $As_3^2:=$ $\left(\begin{array}{ccccc} 1&0&0&0\\
0&0&0&0
\end{array}\right)$

AXIOM 2 $\Longrightarrow$ $\left\{\begin{array}{ccccccccccccc}
\alpha_1 & = & 1  \ \ \  \beta_1 & = & 0 \ \ \ \gamma_1 & = & 1   \ \ \  \delta_1 &  &   \\
\alpha_2 & = & 0  \ \ \ \beta_2 & = & 0 \ \ \ \gamma_2 & = & 0   \ \ \  \delta_2 &  &    \\
\alpha_3 & = & 0 \ \ \  \beta_3 & = & 0 \ \ \ \gamma_3& = & 0\ \ \  \delta_3 &  &   \\
\alpha_4 & = & 0 \  \ \ \beta_3 & = & 0 \ \ \ \gamma_4 & = & 0\ \ \  \delta_4 &  &  \end{array}\right. $

AXIOM 3 $\Longrightarrow$ $\left\{\begin{array}{ccccccccccccc}
\alpha_1 & = & 1  \ \ \  \beta_1 & = & 0 \ \ \ \gamma_1 & = & 1   \ \ \  \delta_1 &  =& 0  \\
\alpha_2 & = & 0  \ \ \ \beta_2 & = & 0 \ \ \ \gamma_2 & = & 0   \ \ \  \delta_2 &  &    \\
\alpha_3 & = & 0 \ \ \  \beta_3 & = & 0 \ \ \ \gamma_3& = & 0\ \ \  \delta_3 & = & 0  \\
\alpha_4 & = & 0 \  \ \ \beta_3 & = & 0 \ \ \ \gamma_4 & = & 0\ \ \  \delta_4 &  &  \end{array}\right. $

AXIOM 4 $\Longrightarrow$ $\left\{\begin{array}{ccccccccccccc}
\alpha_1 & = & 1  \ \ \  \beta_1 & = & 0 \ \ \ \gamma_1 & = & 1   \ \ \  \delta_1 &  =& 0  \\
\alpha_2 & = & 0  \ \ \ \beta_2 & = & 0 \ \ \ \gamma_2 & = & 0   \ \ \  \delta_2 &  &    \\
\alpha_3 & = & 0 \ \ \  \beta_3 & = & 0 \ \ \ \gamma_3& = & 0\ \ \  \delta_3 & = & 0  \\
\alpha_4 & = & 0 \  \ \ \beta_3 & = & 0 \ \ \ \gamma_4 & = & 0\ \ \  \delta_4 & = &0  \end{array}\right. $
%a_1 = 0, a_2 = 0, a_3 = 0, a_4 = 0, b_1 = 1, b_2 = 0, b_3 = 0, b_4 = 0, c_1 = 0, c_2 = 0, c_3 = 0, c_4 = 0, d_2 = 0, d_3 = 0, d_4 = 0

 \textbf{Case} $\delta_2 =0:$

AXIOM 5 $\Longrightarrow$ $\left\{\begin{array}{ccccccccccccc}
\alpha_1 & = & 1  \ \ \  \beta_1 & = & 0 \ \ \ \gamma_1 & = & 1   \ \ \  \delta_1 &  =& 0  \\
\alpha_2 & = & 0  \ \ \ \beta_2 & = & 0 \ \ \ \gamma_2 & = & 0   \ \ \  \delta_2 & = & 0   \\
\alpha_3 & = & 0 \ \ \  \beta_3 & = & 0 \ \ \ \gamma_3& = & 0\ \ \  \delta_3 & = & 0  \\
\alpha_4 & = & 0 \  \ \ \beta_4 & = & 0 \ \ \ \gamma_4 & = & 0\ \ \  \delta_4 & = &0  \end{array}\right. $

$$D_3^2:= \left\{A=\left(\begin{array}{cccc} 1 & 0 & 0 &0\\ 0 & 0 & 0 &0\end{array}\right),\ \ B=\left(\begin{array}{cccc} 1 & 0 & 0 &0\\ 0& 0& 0 & 0\end{array}\right)\right\}.$$

 \textbf{Case} $\delta_2 =1:$

AXIOM 5 $\Longrightarrow$ $\left\{\begin{array}{ccccccccccccc}
\alpha_1 & = & 1  \ \ \  \beta_1 & = & 0 \ \ \ \gamma_1 & = & 1   \ \ \  \delta_1 &  =& 0  \\
\alpha_2 & = & 0  \ \ \ \beta_2 & = & 0 \ \ \ \gamma_2 & = & 0   \ \ \  \delta_2 & = & 1   \\
\alpha_3 & = & 0 \ \ \  \beta_3 & = & 0 \ \ \ \gamma_3& = & 0\ \ \  \delta_3 & = & 0  \\
\alpha_4 & = & 0 \  \ \ \beta_3 & = & 0 \ \ \ \gamma_4 & = & 0\ \ \  \delta_4 & = &0  \end{array}\right. $

$$D_3^3:= \left\{A=\left(\begin{array}{cccc} 1 & 0 & 0 &0\\ 0 & 0 & 0 &0\end{array}\right),\ \ B=\left(\begin{array}{cccc} 1 & 0 & 0 &0\\ 0& 1& 0 & 0\end{array}\right)\right\}.$$
Note that the diassociative algebras $D_3^2$ and $D_3^3$ are not isomorphic since acting by the automorphism group
$$Aut(As_{3}^2)=Aut\left(\begin{pmatrix}
 			1&0&0&0 \\
 			0&0&0&0
 			\end{pmatrix}\right)=\left\{\left(\begin{array}{lll} 1 & 0\\ 0 & t \end{array}\right)| t\neq 0\right\}$$
to the part $B$ of $D_3^2$ 
$$\left(\begin{array}{cccc} 1 & 0 & 0 &0\\ 0& 0& 0 & 0\end{array}\right)=g^{-1}\left(\begin{array}{cccc} 1 & 0 & 0 &0\\ 0& 1& 0 & 0\end{array}\right)g^{\otimes 2}, \ \mbox{where}\ g=\left(\begin{array}{lll} 1 & 0\\ 0 & t \end{array}\right) $$
we get the system of equations which is inconsistent.

Consider $As_3^3:=\left(\begin{array}{ccccc} 1&0&0&0\\
0&1&0&0
\end{array}\right)$

AXIOM 2 $\Longrightarrow$ $\left\{\begin{array}{ccccccccccccc}
\alpha_1 & = & 1  \ \ \  \beta_1 & = & 0 \ \ \ \gamma_1 & = & 1  \ \ \  \delta_1 & = & 0  \\
\alpha_2 & = & 0  \ \ \ \beta_2 & = & 1 \ \ \ \gamma_2 & = & 0   \ \ \  \delta_2 & = &  1  \\
\alpha_3& = & 0 \ \ \  \beta_3 & = & 0 \ \ \ \gamma_3& = & 0\ \ \  \delta_3 & = &  0 \\
\alpha_4 & = & 0 \  \ \ \beta_3 & = & 0 \ \ \ \gamma_4 & = & 0\ \ \  \delta_4 & = & 0 \end{array}\right. $

AXIOM 3,4,5 $\Longrightarrow$ $\left\{\begin{array}{ccccccccccccc}
\alpha_1 & = & 1  \ \ \  \beta_1 & = & 0 \ \ \ \gamma_1 & = & 1  \ \ \  \delta_1 & = & 0  \\
\alpha_2 & = & 0  \ \ \ \beta_2 & = & 1 \ \ \ \gamma_2 & = & 0   \ \ \  \delta_2 & = &  1  \\
\alpha_3 & = & 0 \ \ \  \beta_3 & = & 0 \ \ \ \gamma_3& = & 0\ \ \  \delta_3 & = &  0 \\
\alpha_4 & = & 0 \  \ \ \beta_3 & = & 0 \ \ \ \gamma_4 & = & 0\ \ \  \delta_4 & = & 0 \end{array}\right. $

$$D_3^4:= \left\{A=\left(\begin{array}{cccc} 1 & 0 & 0 &0\\ 0 & 1 & 0 &0\end{array}\right),\ \ B=\left(\begin{array}{cccc} 1 & 0 & 0 &0\\ 0& 1& 0 & 0\end{array}\right)\right\}.$$

Let $A$ to be $As_3^4:=$ $\left(\begin{array}{ccccc} \frac{1}{2}&0&0&0\\
0&0&\frac{1}{2}&0
\end{array}\right)$

AXIOM 2 $\Longrightarrow$ $\left\{\begin{array}{ccccccccccccc}
\alpha_1 & = & \frac{1}{2}  \ \ \  \beta_1 & = & 0 \ \ \ \gamma_1 & = & \frac{1}{2}  \ \ \  \delta_1 &  &   \\
\alpha_2 & = & 0  \ \ \ \beta_2 & = & 0 \ \ \ \gamma_2 & = & 0   \ \ \  \delta_2 &  &    \\
\alpha_3 & = & 0 \ \ \  \beta_3 & = & \frac{1}{2} \ \ \ \gamma_3& = & 0\ \ \  \delta_3 &  &  \\
\alpha_4 & = & 0 \  \ \ \beta_3 & = & 0 \ \ \ \gamma_4 & = & 0\ \ \  \delta_4 &  &  \end{array}\right. $

AXIOM 3 $\Longrightarrow$ $\left\{\begin{array}{ccccccccccccc}
\alpha_1 & = & \frac{1}{2}  \ \ \  \beta_1 & = & 0 \ \ \ \gamma_1 & = & \frac{1}{2}  \ \ \  \delta_1 &  &   \\
\alpha_2 & = & 0  \ \ \ \beta_2 & = & 0 \ \ \ \gamma_2 & = & 0   \ \ \  \delta_2 &  &    \\
\alpha_3 & = & 0 \ \ \  \beta_3 & = & \frac{1}{2} \ \ \ \gamma_3& = & 0\ \ \  \delta_3 &  &  \\
\alpha_4 & = & 0 \  \ \ \beta_3 & = & 0 \ \ \ \gamma_4 & = & 0\ \ \  \delta_4 &  &  \end{array}\right. $

AXIOM 4 gives $\left\{\begin{array}{rrr}\delta_1\delta_3&=&0\\ \delta_1\delta_4&=&0\\ \delta_2\delta_3&=&0\\ \delta_2\delta_4&=&0\\ 2\delta_3^2-\delta_3&=&0\\ 2\delta_3\delta_4-\delta_4&=&0\\ \delta_3\delta_4&=&0\\ \delta_4&=&0.\end{array}\right.$

\textbf{Case 1}: $\delta_3=0$

AXIOM 4 $\Longrightarrow$ $\left\{\begin{array}{ccccccccccccc}
\alpha_1 & = & \frac{1}{2}  \ \ \  \beta_1 & = & 0 \ \ \ \gamma_1 & = & \frac{1}{2}  \ \ \  \delta_1 &  &   \\
\alpha_2 & = & 0  \ \ \ \beta_2 & = & 0 \ \ \ \gamma_2 & = & 0   \ \ \  \delta_2 &  &    \\
\alpha_3 & = & 0 \ \ \  \beta_3 & = & \frac{1}{2} \ \ \ \gamma_3& = & 0\ \ \  \delta_3 & = & 0 \\
\alpha_4 & = & 0 \  \ \ \beta_3 & = & 0 \ \ \ \gamma_4 & = & 0\ \ \  \delta_4 & = & 0 \end{array}\right. $

AXIOM 5 $\Longrightarrow$ $\delta_1 \delta_2=0$ and $\delta_2(2\delta_2-1)=0.$

\textbf{Case 11}: $\delta_2=0$ %a_1 = 1/2, a_2 = 0, a_3 = 0, a_4 = 0, b_1 = 0, b_2 = 0, b_3 = 1/2, b_4 = 0, c_1 = 1/2, c_2 = 0, c_3 = 0, c_4 = 0, d_2 = 0, d_3 = 0, d_4 = 0

AXIOM 5 $\Longrightarrow$ $\left\{\begin{array}{ccccccccccccc}
\alpha_1 & = & \frac{1}{2}  \ \ \  \beta_1 & = & 0 \ \ \ \gamma_1 & = & \frac{1}{2}  \ \ \  \delta_1 &  &   \\
\alpha_2 & = & 0  \ \ \ \beta_2 & = & 0 \ \ \ \gamma_2 & = & 0   \ \ \  \delta_2 & = &  0  \\
\alpha_3 & = & 0 \ \ \  \beta_3 & = & \frac{1}{2} \ \ \ \gamma_3& = & 0\ \ \  \delta_3 & = & 0 \\
\alpha_4 & = & 0 \  \ \ \beta_3 & = & 0 \ \ \ \gamma_4 & = & 0\ \ \  \delta_4 & = & 0 \end{array}\right. $

$$D_3^5:= \left\{A=\left(\begin{array}{cccc} \frac{1}{2} & 0 & 0 &0\\ 0 & 0 & \frac{1}{2} &0\end{array}\right),\ \ B=\left(\begin{array}{cccc} \frac{1}{2} & 0 & 0 &0\\ \delta_1& 0& 0 & 0\end{array}\right)\right\}.$$

\textbf{Case 12}: \ \ $\delta_2=\frac{1}{2},\ \delta_1=0$

AXIOM 5 $\Longrightarrow$ $\left\{\begin{array}{ccccccccccccc}
\alpha_1 & = & \frac{1}{2}  \ \ \  \beta_1 & = & 0 \ \ \ \gamma_1 & = & \frac{1}{2}  \ \ \  \delta_1 & = & 0  \\
\alpha_2 & = & 0  \ \ \ \beta_2 & = & 0 \ \ \ \gamma_2 & = & 0   \ \ \  \delta_2 & = &  \frac{1}{2}  \\
\alpha_3 & = & 0 \ \ \  \beta_3 & = & \frac{1}{2} \ \ \ \gamma_3& = & 0\ \ \  \delta_3 & = & 0 \\
\alpha_4 & = & 0 \  \ \ \beta_3 & = & 0 \ \ \ \gamma_4 & = & 0\ \ \  \delta_4 & = & 0 \end{array}\right. $

$$D_3^6:= \left\{A=\left(\begin{array}{cccc} \frac{1}{2} & 0 & 0 &0\\ 0 & 0 & \frac{1}{2} &0\end{array}\right),\ \ B=\left(\begin{array}{cccc} \frac{1}{2} & 0 & 0 &0\\ 0&\frac{1}{2} & 0 & 0\end{array}\right)\right\}.$$

The dialgebras $D_3^5$ and $D_3^6$ are not isomorphic that can be seen acting by $$Aut(As_{3}^4)=Aut\left(\begin{pmatrix}
 			\frac{1}{2}&0&0&0 \\
 			0&0&\frac{1}{2}&0
 			\end{pmatrix}\right)=\left\{\left(\begin{array}{lll} 1 & 0\\ z & t \end{array}\right)\big| \ t\neq 0\right\}$$
 
 on the part $B$ of $D_3^5$
 $$\left(\begin{array}{cccc} \frac{1}{2} & 0 & 0 &0\\ \delta_1& 0& 0 & 0\end{array}\right)=g^{-1}\left(\begin{array}{cccc} \frac{1}{2} & 0 & 0 &0\\ 0& \frac{1}{2}& 0 & 0\end{array}\right)g^{\otimes 2}, \ \mbox{where}\ g=\left(\begin{array}{lll} 1 & 0\\ z & t \end{array}\right). $$
We get an inconsistent system of equations.

\textbf{Case 2}: \ \ $\delta_3\neq 0$ $\Longrightarrow$ $\delta_1=0$ \ $\delta_2=0$ and $\delta_3=\frac{1}{2}$

AXIOM 4,5 $\Longrightarrow$ $\left\{\begin{array}{ccccccccccccc}
\alpha_1 & = & \frac{1}{2}  \ \ \  \beta_1 & = & 0 \ \ \ \gamma_1 & = & \frac{1}{2}  \ \ \  \delta_1 & = & 0  \\
\alpha_2 & = & 0  \ \ \ \beta_2 & = & 0 \ \ \ \gamma_2 & = & 0   \ \ \  \delta_2 & = &  0  \\
\alpha_3 & = & 0 \ \ \  \beta_3 & = & \frac{1}{2} \ \ \ \gamma_3& = & 0\ \ \  \delta_3 & = & \frac{1}{2} \\
\alpha_4 & = & 0 \  \ \ \beta_3 & = & 0 \ \ \ \gamma_4 & = & 0\ \ \  \delta_4 & = & 0 \end{array}\right. $

$$D_3^7:= \left\{A=\left(\begin{array}{cccc} \frac{1}{2} & 0 & 0 &0\\ 0 & 0 & \frac{1}{2} &0\end{array}\right),\ \ B=\left(\begin{array}{cccc} \frac{1}{2} & 0 & 0 &0\\ 0&0 & \frac{1}{2} & 0\end{array}\right)\right\}.$$

Consider $A_3^5(\alpha_4):=$ $\left(\begin{array}{ccccc} \frac{1}{2}&0&0&\alpha_4\\
0&\frac{1}{2}&\frac{1}{2}&0
\end{array}\right)$, $\alpha_4 \in \mathbb{F}.$

AXIOM 2 $\Longrightarrow$ $\left\{\begin{array}{ccccccccccccc}
\alpha_1 & = & \frac{1}{2}  \ \ \  \beta_1 & = & 0 \ \ \ \gamma_1 & = & \frac{1}{2}  \ \ \  \delta_1 & = & 0  \\
\alpha_2 & = & 0  \ \ \ \beta_2 & = & \frac{1}{2} \ \ \ \gamma_2 & = & 0   \ \ \  \delta_2 & = &  \frac{1}{2}  \\
\alpha_3 & = & 0 \ \ \  \beta_3 & = & \frac{1}{2} \ \ \ \gamma_3& = & 0\ \ \  \delta_3 & = & \frac{1}{2} \\
\alpha_4 &  &  \  \ \ \beta_3 & = & 0 \ \ \ \gamma_4 & = & \alpha_4\ \ \  \delta_4 & = & 0 \end{array}\right. $

AXIOM 3,4,5 $\Longrightarrow$ $\left\{\begin{array}{ccccccccccccc}
\alpha_1 & = & \frac{1}{2}  \ \ \  \beta_1 & = & 0 \ \ \ \gamma_1 & = & \frac{1}{2}  \ \ \  \delta_1 & = & 0  \\
\alpha_2 & = & 0  \ \ \ \beta_2 & = & \frac{1}{2} \ \ \ \gamma_2 & = & 0   \ \ \  \delta_2 & = &  \frac{1}{2}  \\
\alpha_3 & = & 0 \ \ \  \beta_3 & = & \frac{1}{2} \ \ \ \gamma_3& = & 0\ \ \  \delta_3 & = & \frac{1}{2} \\
\alpha_4 &  &  \  \ \ \beta_3 & = & 0 \ \ \ \gamma_4 & = & \alpha_4\ \ \  \delta_4 & = & 0 \end{array}\right. $

$$D_3^8:= \left\{A=\left(\begin{array}{cccc} \frac{1}{2} & 0 & 0 &\alpha_4\\ 0 & \frac{1}{2} & \frac{1}{2} &0\end{array}\right),\ \ B=\left(\begin{array}{cccc} \frac{1}{2} & 0 & 0 &\alpha_4\\ 0&\frac{1}{2} & \frac{1}{2} & 0\end{array}\right)\right\}.$$

\begin{theor} \label{char0Dias} Any non-trivial $2$-dimensional associative dialgebra over a field $\mathbb{F},$ $(Char(\mathbb{F})\neq 2,3)$ is isomorphic to only one of the following listed by their matrices of structure constants, such algebras:

\begin{enumerate}
\item Diassociative algebras generated by $A_{13}$:
%a_1 = 0, a_2 = 0, a_3 = 0, a_4 = 0, b_1 = 1, b_2 = 0, b_3 = 0, b_4 = 0, c_1 = 0, c_2 = 0, c_3 = 0, c_4 = 0, d_2 = 0, d_3 = 0, d_4 = 0
\begin{itemize}
\item $D_{13}^1:=\left\{A=\left(\begin{array}{cccc} 0 & 0 & 0 &0\\ 1 & 0 & 0 &0\end{array}\right),\ \ B=\left(\begin{array}{cccc} 0 & 0 & 0 &0\\ \delta_1& 0& 0 & 0\end{array}\right)\right\}$
\end{itemize}
\item Diassociative algebras generated by $A_{3}$:

\begin{itemize}%[a_1 = 1, a_2 = 0, a_3 = 0, a_4 = 0, b_1 = 0, b_2 = 0, b_3 = 0, b_4 = 0, c_1 = 1, c_2 = 0, c_3 = 0, c_4 = 0, d_1 = 0, d_2 = 0, d_3 = 0, d_4 = 0]
\item $D_3^2:=\left\{A=\left(\begin{array}{cccc} 1 & 0 & 0 &0\\ 0 & 0 & 0 &0\end{array}\right),\ \ B=\left(\begin{array}{cccc} 1 & 0 & 0 &0\\ 0& 0& 0 & 0\end{array}\right)\right\}$
\item %a_1 = 1, a_2 = 0, a_3 = 0, a_4 = 0, b_1 = 0, b_2 = 0, b_3 = 0, b_4 = 0, c_1 = 1, c_2 = 0, c_3 = 0, c_4 = 0, d_1 = 0, d_2 = 1, d_3 = 0, d_4 = 0
$D_3^3:=\left\{A=\left(\begin{array}{cccc} 1 & 0 & 0 &0\\ 0 & 0 & 0 &0\end{array}\right),\ \ B=\left(\begin{array}{cccc} 1 & 0 & 0 &0\\ 0& 1& 0 & 0\end{array}\right)\right\}$
\item %a_1 = 1, a_2 = 0, a_3 = 0, a_4 = 0, b_1 = 0, b_2 = 1, b_3 = 0, b_4 = 0, c_1 = 1, c_2 = 0, c_3 = 0, c_4 = 0, d_1 = 0, d_2 = 1, d_3 = 0, d_4 = 0
$D_3^4:=\left\{A=\left(\begin{array}{cccc} 1 & 0 & 0 &0\\ 0 & 1 & 0 &0\end{array}\right),\ \ B=\left(\begin{array}{cccc} 1 & 0 & 0 &0\\ 0& 1& 0 & 0\end{array}\right)\right\}$
\item %a_1 = 1/2, a_2 = 0, a_3 = 0, a_4 = 0, b_1 = 0, b_2 = 0, b_3 = 1/2, b_4 = 0, c_1 = 1/2, c_2 = 0, c_3 = 0, c_4 = 0, d_2 = 0, d_3 = 0, d_4 = 0
$D_3^5(\delta_1):=\left\{A=\left(\begin{array}{cccc} \frac{1}{2} & 0 & 0 &0\\ 0 & 0 & \frac{1}{2} &0\end{array}\right),\ \ B=\left(\begin{array}{cccc} \frac{1}{2} & 0 & 0 &0\\ \delta_1& 0& 0 & 0\end{array}\right), \ \ \delta_1\in \mathbb{F}\right\}$
\item %a_1 = 1/2, a_2 = 0, a_3 = 0, a_4 = 0, b_1 = 0, b_2 = 0, b_3 = 1/2, b_4 = 0, c_1 = 1/2, c_2 = 0, c_3 = 0, c_4 = 0, d_1 = 0, d_2 = 1/2, d_3 = 0, d_4 = 0
 $D_3^6:=\left\{A=\left(\begin{array}{cccc} \frac{1}{2} & 0 & 0 &0\\ 0 & 0 & \frac{1}{2} &0\end{array}\right),\ \ B=\left(\begin{array}{cccc} \frac{1}{2} & 0 & 0 &0\\ 0& \frac{1}{2}& 0 & 0\end{array}\right)\right\}$

\item %a_1 = 1/2, a_2 = 0, a_3 = 0, a_4 = 0, b_1 = 0, b_2 = 0, b_3 = 1/2, b_4 = 0, c_1 = 1/2, c_2 = 0, c_3 = 0, c_4 = 0, d_1 = 0, d_2 = 0, d_3 = 1/2, d_4 = 0
$D_3^7:=\left\{A=\left(\begin{array}{cccc} \frac{1}{2} & 0 & 0 &0\\ 0 & 0 & \frac{1}{2} &0\end{array}\right),\ \ B=\left(\begin{array}{cccc} \frac{1}{2} & 0 & 0 &0\\ 0& 0& \frac{1}{2} & 0\end{array}\right)\right\}$

\item %a_1 = 1/2, a_2 = 0, a_3 = 0, b_1 = 0, b_2 = 1/2, b_3 = 1/2, b_4 = 0, c_1 = 1/2, c_2 = 0, c_3 = 0, c_4 = a_4, d_1 = 0, d_2 = 1/2, d_3 = 1/2, d_4 = 0
$D_3^8:=\left\{A=\left(\begin{array}{cccc} \frac{1}{2} & 0 & 0 &\alpha_4\\ 0 & \frac{1}{2} & \frac{1}{2} &0\end{array}\right),\ \ B=\left(\begin{array}{cccc} \frac{1}{2} & 0 & 0 &\alpha_4\\ 0& \frac{1}{2}& \frac{1}{2} & 0\end{array}\right), \ \ \alpha_4 \in \mathbb{F}\right\}$
\end{itemize}
\end{enumerate}
\end{theor}

According to a result of \cite{WI} there are four classes of two-dimensional associative dialgebras over $\mathbb{C}$ given as follows

$Dias^1:=\left\{A=\left(\begin{array}{cccc}1 & 0 & 0 &0\\ 0& 0& 1 & 0\end{array}\right),\ \ B=\left(\begin{array}{cccc}1 & 0 & 0 &0\\ 0 & 0 & 0 &0\end{array}\right)\right\}\cong D_3^5(0);$

$Dias^2:=\left\{ A=\left(\begin{array}{cccc}1 & 0 & 0 &0\\ 0& 0& 0 & 0\end{array}\right),\ \ B=\left(\begin{array}{cccc}1 & 0 & 0 &0\\ 0 & 1 & 0 &0\end{array}\right)\right\}\cong D_3^3;$

$Dias^3:=\left\{A=\left(\begin{array}{cccc}0 & 0 & 0 &0\\ 1 & 0 & 0 &0\end{array}\right),\ \ B=\left(\begin{array}{cccc}0 & 0 & 0 &0\\ \alpha& 0& 0 & 0\end{array}\right),\ \ \alpha\in \mathbb{C}\right\}\cong D_{13}^1;$

$Dias^4:=\left\{ A=\left(\begin{array}{cccc}1 & 0 & 0 &0\\ 0& 0& 1 & 0\end{array}\right),\ \ B=\left(\begin{array}{cccc}1 & 0 & 0 &0\\ 0 & 1 & 0 &0\end{array}\right)\right\}\cong D_3^6.$

Since Theorem \ref{char0Dias} includes the case $Char(\mathbb{C})=0$  the list in \cite{WI} must be accordingly corrected.

\subsection{Characteristic of $\mathbb{F}$ is two}\emph{}

In the case of the characteristic of the field $\mathbb{F}$ is two the associative dialgebras generated from the list of Theorem \ref{char2Asso} are as follows:

From $As_{12,2}^1=\begin{pmatrix}
 		0&0&0&0 \\
 		1&0&0&0\\
 		\end{pmatrix} 
 		$  we get

%A=[a_1 = 0, a_2 = 0, a_3 = 0, a_4 = 0, b_1 = 1, b_2 = 0, b_3 = 0, b_4 = 0]
%AB=[a_1 = 0, a_2 = 0, a_3 = 0, a_4 = 0, b_1 = 1, b_2 = 0, b_3 = 0, b_4 = 0, c_1 = 0, c_2 = 0, c_3 = 0, c_4 = 0, d_2 = 0, d_3 = 0, d_4 = 0]

$$D_{12,2}^1:=\left\{A=\left(\begin{array}{cccc}  0& 0 & 0 &0\\ 1 & 0 & 0 &0\end{array}\right),\ \ B=\left(\begin{array}{cccc} 0 & 0 & 0 &0\\ \delta_1& 0& 0 & 0\end{array}\right), \ \ \delta_1 \in \mathbb{F}\right\}.$$

The algebra $As_{11,2}^2(\beta_1)=\begin{pmatrix}
 		0&1&1&0 \\
 		\beta_1&0&0&1\\
 		\end{pmatrix}\simeq\begin{pmatrix}
 		0&1&1&0 \\
 		b^2(\beta_1+a^2)&0&0&1\\
 		\end{pmatrix}
 		,$ where $a, b, \beta_1\in\mathbb{F}$ and $ b\neq 0$ produces

%A=[a_1 = 0, a_2 = 1, a_3 = 1, a_4 = 0, b_2 = 0, b_3 = 0, b_4 = 1]
%AB=A11_5_d20d30c10d41c21c31c40d1b1 := subs([a_1 = 0, a_2 = 1, a_3 = 1, a_4 = 0, b_2 = 0, b_3 = 0, b_4 = 1, c_1 = 0, c_2 = 1, c_3 = 1, c_4 = 0, d_1 = b_1, d_2 = 0, d_3 = 0, d_4 = 1], AX5)

$$D_{11,2}^2:=\left\{A=\left(\begin{array}{cccc} 0 & 1 & 1 &0\\ \beta_1 & 0 & 0 &1\end{array}\right),\ \ B=\left(\begin{array}{cccc} 0 & 1 & 1 &0\\ \beta_1& 0& 0 & 1\end{array}\right), \ \ \beta_1 \in \mathbb{F}\right\}.$$

From $As_{6,2}^3=\begin{pmatrix}
 		1&0&0&0 \\
 		0& 0& 1& 0
 		\end{pmatrix}$ we get

%A=A6_1_a1n_a40 := subs([a_1 = -1, a_2 = 0, a_3 = 0, a_4 = 0, b_1 = 0, b_2 = 0, b_3 = -1, b_4 = 0], AX1)
%AB=A6_5_a1n_a40_c_c1n1_d40_d30_d20 := subs([a_1 = -1, a_2 = 0, a_3 = 0, a_4 = 0, b_1 = 0, b_2 = 0, b_3 = -1, b_4 = 0, c_1 = -1, c_2 = 0, c_3 = 0, c_4 = 0, d_2 = 0, d_3 = 0, d_4 = 0], AX5)
$$D_{6,2}^3:=\left\{A=\left(\begin{array}{cccc} 1 & 0 & 0 &0\\ 0 & 0 & 1 &0\end{array}\right),\ \ B=\left(\begin{array}{cccc} 1 & 0 & 0 &0\\ \beta_1& 0& 0 & 0\end{array}\right), \ \ \beta_1 \in \mathbb{F}\right\}$$
%AB=A6_5_a1n_a40_c_c1n1_d40_d30_d10d21 := subs([a_1 = -1, a_2 = 0, a_3 = 0, a_4 = 0, b_1 = 0, b_2 = 0, b_3 = -1, b_4 = 0, c_1 = -1, c_2 = 0, c_3 = 0, c_4 = 0, d_1 = 0, d_2 = 1, d_3 = 0, d_4 = 0], AX5)
$$D_{6,2}^4:=\left\{A=\left(\begin{array}{cccc} 1 & 0 & 0 &0\\ 0 & 0 & 1 &0\end{array}\right),\ \ B=\left(\begin{array}{cccc} 1 & 0 & 0 &0\\ 0& 1& 0 & 0\end{array}\right) \right\}$$
%AB=A6_5_a1n_a40_c_c1n1_d40 := subs([a_1 = -1, a_2 = 0, a_3 = 0, a_4 = 0, b_1 = 0, b_2 = 0, b_3 = -1, b_4 = 0, c_1 = -1, c_2 = 0, c_3 = 0, c_4 = 0, d_1 = 0, d_2 = 0, d_3 = 1, d_4 = 0], AX5)
$$D_{6,2}^5:=\left\{A=\left(\begin{array}{cccc} 1 & 0 & 0 &0\\ 0 & 0 & 1 &0\end{array}\right),\ \ B=\left(\begin{array}{cccc} 1 & 0 & 0 &0\\ 0& 0& 1 & 0\end{array}\right) \right\}.$$

The diassociative algebras $D_{6,2}^3$, $D_{6,2}^4$ and $D_{6,2}^5$ are not isomorphic to each others since the group of automorphisms of $As_{6,2}^3$ is trivial.

Consider $As_{4,2}^4(\beta_1)=\begin{pmatrix}
 		1&1&1&0 \\
 		\beta_1& 0& 0& 1
 		\end{pmatrix}.$ This generates
%A4_b20_a1n1 := subs([a_1 = -1, a_2 = 1, a_3 = 1, a_4 = 0, b_2 = 0, b_3 = 0, b_4 = 1], AX)
%A4_5_b20_a1n1_c11c21c31c4n_b10 := subs([a_1 = -1, a_2 = 1, a_3 = 1, a_4 = 0, b_2 = 0, b_1 = 0, b_3 = 0, b_4 = 1, c_1 = -1, c_2 = -1, c_3 = 1, c_4 = 0, d_1 = 0, d_2 = 0, d_3 = 0, d_4 = 1], AX5)
$$D_{4,2}^6:=\left\{A=\left(\begin{array}{cccc} 1 & 1 & 1 &0\\ 0 & 0 & 0 &1\end{array}\right),\ \ B=\left(\begin{array}{cccc} 1 & 1 & 1 &0\\ 0& 0& 0 & 1\end{array}\right) \right\}.$$

From $ As_{3,2}^5=\begin{pmatrix}
 		1&0&0&0 \\
 		0& 0& 0& 0
 		\end{pmatrix}$ we get

%A=A3_1_a11_a40_b20 := subs([a_1 = 1, a_2 = 0, a_3 = 0, a_4 = 0, b_1 = 0, b_2 = 0, b_3 = 0, b_4 = 0], AX1)
%AB=A3_5_a11_a40_b20_c_c11_d10d30_d40_d20 := subs([a_1 = 1, a_2 = 0, a_3 = 0, a_4 = 0, b_1 = 0, b_2 = 0, b_3 = 0, b_4 = 0, c_1 = 1, c_2 = 0, c_3 = 0, c_4 = 0, d_1 = 0, d_2 = 0, d_3 = 0, d_4 = 0], AX5)
$$D_{3,2}^7:=\left\{A=\left(\begin{array}{cccc} 1 & 0 & 0 &0\\ 0 & 0 & 0 &0\end{array}\right),\ \ B=\left(\begin{array}{cccc} 1 & 0 & 0 &0\\ 0& 0& 0 & 0\end{array}\right) \right\}$$

%AB=A3_5_a11_a40_b20_c_c11_d10d30_d40_d21 := subs([a_1 = 1, a_2 = 0, a_3 = 0, a_4 = 0, b_1 = 0, b_2 = 0, b_3 = 0, b_4 = 0, c_1 = 1, c_2 = 0, c_3 = 0, c_4 = 0, d_1 = 0, d_2 = 1, d_3 = 0, d_4 = 0], AX5)
$$D_{3,2}^8:=\left\{A=\left(\begin{array}{cccc} 1 & 0 & 0 &0\\ 0 & 0 & 0 &0\end{array}\right),\ \ B=\left(\begin{array}{cccc} 1 & 0 & 0 &0\\ 0& 1& 0 & 0\end{array}\right) \right\}$$

The algebras $D_{3,2}^7$ and $D_{3,2}^8$ are not isomorphic since there is no an element of the automorphism group $$Aut(As_{3,2}^6)=\left\{ \left(\begin{array}{crr} 1 & 0\\ z & t \end{array}\right)\big| \ \ z,t \in \mathbb{F} \ \mbox{and}\ t \neq 0 \right\}.$$
sending the part $B$ of $D_{3,2}^7$ to the part $B$ of $D_{3,2}^8$.

Finally, from $ As_{3,2}^6=\begin{pmatrix}
 		1&0&0&0 \\
 		0& 1& 0& 0
 		\end{pmatrix}$ we get
%A=A3_1_a11_a40_b21 := subs([a_1 = 1, a_2 = 0, a_3 = 0, a_4 = 0, b_1 = 0, b_2 = 1, b_3 = 0, b_4 = 0], AX1)
%%AB=A3_5_a11_a40_b21_c_c11_d_d21 := subs([a_1 = 1, a_2 = 0, a_3 = 0, a_4 = 0, b_1 = 0, b_2 = 1, b_3 = 0, b_4 = 0, c_1 = 1, c_2 = 0, c_3 = 0, c_4 = 0, d_1 = 0, d_2 = 1, d_3 = 0, d_4 = 0], AX5)
$$D_{3,2}^9:=\left\{A=\left(\begin{array}{cccc} 1 & 0 & 0 &0\\ 0 & 1 & 0 &0\end{array}\right),\ \ B=\left(\begin{array}{cccc} 1 & 0 & 0 &0\\ 0& 1& 0 & 0\end{array}\right) \right\}.$$
\begin{theor} \label{char2Dias} Any non-trivial $2$-dimensional associative dialgebra over a field $\mathbb{F},$ $(Char(\mathbb{F})= 2)$ is isomorphic to only one of the following listed by their matrices of structure constants, such algebras:

\begin{enumerate}
%\item Diassociative algebras generated by $A_{12,2}^1$:
%a_1 = 0, a_2 = 0, a_3 = 0, a_4 = 0, b_1 = 1, b_2 = 0, b_3 = 0, b_4 = 0, c_1 = 0, c_2 = 0, c_3 = 0, c_4 = 0, d_2 = 0, d_3 = 0, d_4 = 0
%\begin{itemize}
\item $D_{12,2}^1:=\left\{A=\left(\begin{array}{cccc}  0& 0 & 0 &0\\ 1 & 0 & 0 &0\end{array}\right),\ \ B=\left(\begin{array}{cccc} 0 & 0 & 0 &0\\ \delta_1& 0& 0 & 0\end{array}\right), \ \ \delta_1 \in \mathbb{F}\right\}$
%\end{itemize}
%
%\item Diassociative algebras generated by $A_{3}$:

%\begin{itemize}%[a_1 = 1, a_2 = 0, a_3 = 0, a_4 = 0, b_1 = 0, b_2 = 0, b_3 = 0, b_4 = 0, c_1 = 1, c_2 = 0, c_3 = 0, c_4 = 0, d_1 = 0, d_2 = 0, d_3 = 0, d_4 = 0]
\item $D_{11,2}^2:=\left\{A=\left(\begin{array}{cccc} 0 & 1 & 1 &0\\ \beta_1 & 0 & 0 &1\end{array}\right),\ \ B=\left(\begin{array}{cccc} 0 & 1 & 1 &0\\ \beta_1& 0& 0 & 1\end{array}\right), \ \ \beta_1 \in \mathbb{F}\right\}$
    
    %$$D_{6,2}^3:=\left\{A=\left(\begin{array}{cccc} 1 & 0 & 0 &0\\ 0 & 0 & 1 &0\end{array}\right)\ \ B=\left(\begin{array}{cccc} 1 & 0 & 0 &0\\ \beta_1& 0& 0 & 0\end{array}\right), \ \ \beta_1 \in \mathbb{F}\right\}$$
%AB=A6_5_a1n_a40_c_c1n1_d40_d30_d10d21 := subs([a_1 = -1, a_2 = 0, a_3 = 0, a_4 = 0, b_1 = 0, b_2 = 0, b_3 = -1, b_4 = 0, c_1 = -1, c_2 = 0, c_3 = 0, c_4 = 0, d_1 = 0, d_2 = 1, d_3 = 0, d_4 = 0], AX5)
%$$$$
%AB=A6_5_a1n_a40_c_c1n1_d40 := subs([a_1 = -1, a_2 = 0, a_3 = 0, a_4 = 0, b_1 = 0, b_2 = 0, b_3 = -1, b_4 = 0, c_1 = -1, c_2 = 0, c_3 = 0, c_4 = 0, d_1 = 0, d_2 = 0, d_3 = 1, d_4 = 0], AX5)
%$$$$
\item %a_1 = 1, a_2 = 0, a_3 = 0, a_4 = 0, b_1 = 0, b_2 = 0, b_3 = 0, b_4 = 0, c_1 = 1, c_2 = 0, c_3 = 0, c_4 = 0, d_1 = 0, d_2 = 1, d_3 = 0, d_4 = 0
$D_{6,2}^3:=\left\{A=\left(\begin{array}{cccc} 1 & 0 & 0 &0\\ 0 & 0 & 1 &0\end{array}\right),\ \ B=\left(\begin{array}{cccc} 1 & 0 & 0 &0\\ \beta_1& 0& 0 & 0\end{array}\right), \ \ \beta_1 \in \mathbb{F}\right\}$
\item %a_1 = 1, a_2 = 0, a_3 = 0, a_4 = 0, b_1 = 0, b_2 = 1, b_3 = 0, b_4 = 0, c_1 = 1, c_2 = 0, c_3 = 0, c_4 = 0, d_1 = 0, d_2 = 1, d_3 = 0, d_4 = 0
$D_{6,2}^4:=\left\{A=\left(\begin{array}{cccc} 1 & 0 & 0 &0\\ 0 & 0 & 1 &0\end{array}\right),\ \ B=\left(\begin{array}{cccc} 1 & 0 & 0 &0\\ 0& 1& 0 & 0\end{array}\right) \right\}$
\item %a_1 = 1/2, a_2 = 0, a_3 = 0, a_4 = 0, b_1 = 0, b_2 = 0, b_3 = 1/2, b_4 = 0, c_1 = 1/2, c_2 = 0, c_3 = 0, c_4 = 0, d_2 = 0, d_3 = 0, d_4 = 0
$D_{6,2}^5:=\left\{A=\left(\begin{array}{cccc} 1 & 0 & 0 &0\\ 0 & 0 & 1 &0\end{array}\right),\ \ B=\left(\begin{array}{cccc} 1 & 0 & 0 &0\\ 0& 0& 1 & 0\end{array}\right) \right\}$
\item %a_1 = 1/2, a_2 = 0, a_3 = 0, a_4 = 0, b_1 = 0, b_2 = 0, b_3 = 1/2, b_4 = 0, c_1 = 1/2, c_2 = 0, c_3 = 0, c_4 = 0, d_1 = 0, d_2 = 1/2, d_3 = 0, d_4 = 0
 $D_{4,2}^6:=\left\{A=\left(\begin{array}{cccc} 1 & 1 & 1 &0\\ 0 & 0 & 0 &1\end{array}\right),\ \ B=\left(\begin{array}{cccc} 1 & 1 & 1 &0\\ 0& 0& 0 & 1\end{array}\right) \right\}$

%$$$$

%AB=A3_5_a11_a40_b20_c_c11_d10d30_d40_d21 := subs([a_1 = 1, a_2 = 0, a_3 = 0, a_4 = 0, b_1 = 0, b_2 = 0, b_3 = 0, b_4 = 0, c_1 = 1, c_2 = 0, c_3 = 0, c_4 = 0, d_1 = 0, d_2 = 1, d_3 = 0, d_4 = 0], AX5)
%$$$$

\item %a_1 = 1/2, a_2 = 0, a_3 = 0, a_4 = 0, b_1 = 0, b_2 = 0, b_3 = 1/2, b_4 = 0, c_1 = 1/2, c_2 = 0, c_3 = 0, c_4 = 0, d_1 = 0, d_2 = 0, d_3 = 1/2, d_4 = 0
$D_{3,2}^7:=\left\{A=\left(\begin{array}{cccc} 1 & 0 & 0 &0\\ 0 & 0 & 0 &0\end{array}\right),\ \ B=\left(\begin{array}{cccc} 1 & 0 & 0 &0\\ 0& 0& 0 & 0\end{array}\right) \right\}$

\item %a_1 = 1/2, a_2 = 0, a_3 = 0, b_1 = 0, b_2 = 1/2, b_3 = 1/2, b_4 = 0, c_1 = 1/2, c_2 = 0, c_3 = 0, c_4 = a_4, d_1 = 0, d_2 = 1/2, d_3 = 1/2, d_4 = 0
$D_{3,2}^8:=\left\{A=\left(\begin{array}{cccc} 1 & 0 & 0 &0\\ 0 & 0 & 0 &0\end{array}\right),\ \ B=\left(\begin{array}{cccc} 1 & 0 & 0 &0\\ 0& 1& 0 & 0\end{array}\right) \right\}$
%\end{itemize}
\item $D_{3,2}^9:=\left\{A=\left(\begin{array}{cccc} 1 & 0 & 0 &0\\ 0 & 1 & 0 &0\end{array}\right),\ \ B=\left(\begin{array}{cccc} 1 & 0 & 0 &0\\ 0& 1& 0 & 0\end{array}\right) \right\}.$
\end{enumerate}
\end{theor}

\subsection{Characteristic of $\mathbb{F}$ is three}\emph{}

The associative algebra $As_{13,3}^1=\begin{pmatrix}
 		0&0&0&0 \\
 		1&0&0&0
 		\end{pmatrix}$ produces the diassociative algebra

%A=A13_1 := subs([a_1 = 0, a_2 = 0, a_3 = 0, a_4 = 0, b_1 = 1, b_2 = 0, b_3 = 0, b_4 = 0], AX1);
%AB=A13_5_c_d20d40_d30 := subs([a_1 = 0, a_2 = 0, a_3 = 0, a_4 = 0, b_1 = 1, b_2 = 0, b_3 = 0, b_4 = 0, c_1 = 0, c_2 = 0, c_3 = 0, c_4 = 0, d_2 = 0, d_3 = 0, d_4 = 0], AX5)
$$D_{13,3}^1:=\left\{A=\left(\begin{array}{cccc} 0 & 0 & 0 &0\\ 1 & 0 & 0 &0\end{array}\right),\ \ B=\left(\begin{array}{cccc} 0 & 0 & 0 &0\\ \delta_1& 0& 0 & 0\end{array}\right), \ \ \delta_1 \in \mathbb{F} \right\}.$$

%From $A_{11,3}(\mathrm{c})=\begin{pmatrix}
% 		0&1&1&0 \\
% 		\beta_1&0&0&-1\\
% 		\end{pmatrix}\simeq\begin{pmatrix}
% 		0&1&1&0 \\
% 		a^{2}\beta_1&0&0&-1\\
% 		\end{pmatrix},
% 		$ where  $a, \mathrm{c}=\beta_1\in\mathbb{F},$ $a\neq 0$
%
%%A=A11_1 := subs([a_1 = 0, a_2 = 1, a_3 = 1, a_4 = 0, b_2 = 0, b_3 = 0, b_4 = 1], AX1)
%%AB=A11_5_c10c21c31c40d1b1d20d30d41 := subs([a_1 = 0, a_2 = 1, a_3 = 1, a_4 = 0, b_2 = 0, b_3 = 0, b_4 = 1, c_1 = 0, c_2 = 1, c_3 = 1, c_4 = 0, d_1 = b_1, d_2 = 0, d_3 = 0, d_4 = 1], AX5)
%$$D_{11,3}^2:=\left\{A=\left(\begin{array}{cccc} 0 & 1 & 1 &0\\ \beta_1 & 0 & 0 &1\end{array}\right)\ \ B=\left(\begin{array}{cccc} 0 & 1 & 1 &0\\ \beta_1& 0& 0 & 1\end{array}\right), \ \ \beta_1 \in \mathbb{F} \right\}$$
%
%From $A_{10,3}(\mathrm{c})=\begin{pmatrix}
% 		0&0&0&1\\
% 		\beta_1&0&0&0
% 		\end{pmatrix}\simeq\begin{pmatrix}
% 		0&0&0&1\\
% 		a^{3}\beta_1^{\pm1}&0&0&0
% 		\end{pmatrix},$\\  where the polynomial $ \beta_1-t^3$ has no root, $ a, \mathrm{c}=\beta_1\in\mathbb{F}$ and  $ a, \beta_1\neq 0$
%
%%A=A10_1 := subs([a_1 = 0, a_2 = 0, a_3 = 0, a_4 = 0, b_2 = 0, b_3 = 0, b_4 = 0], AX1)
%
%%AB=A10_5_c_d20d40_d30 := subs([a_1 = 0, a_2 = 0, a_3 = 0, a_4 = 0, b_2 = 0, b_3 = 0, b_4 = 0, c_1 = 0, c_2 = 0, c_3 = 0, c_4 = 0, d_2 = 0, d_3 = 0, d_4 = 0], AX5)
%$$D_{10,3}^3:=\left\{A=\left(\begin{array}{cccc} 0 & 0 & 0 &0\\ \beta_1 & 0 & 0 &0\end{array}\right)\ \ B=\left(\begin{array}{cccc} 0 & 0 & 0 &0\\ \delta_1& 0& 0 & 0\end{array}\right), \ \ \beta_1,\delta_1 \in \mathbb{F} \right\}$$

From $As_{3,3}^2=\begin{pmatrix}
 		1&0&0&0 \\
 		0& 0& 0& 0
 		\end{pmatrix}$ we get the diassociative algebras
%A=A3_1_a11a40b20 := subs([a_1 = 1, a_2 = 0, a_3 = 0, a_4 = 0, b_1 = 0, b_2 = 0, b_3 = 0, b_4 = 0], AX1)
%AB=A3_5_a11a40b20_c_c11_d10d30_d40_d20 := subs([a_1 = 1, a_2 = 0, a_3 = 0, a_4 = 0, b_1 = 0, b_2 = 0, b_3 = 0, b_4 = 0, c_1 = 1, c_2 = 0, c_3 = 0, c_4 = 0, d_1 = 0, d_2 = 0, d_3 = 0, d_4 = 0], AX5)

$$D_{3,3}^2:=\left\{A=\left(\begin{array}{cccc} 1 & 0 & 0 &0\\ 0 & 0 & 0 &0\end{array}\right),\ \ B=\left(\begin{array}{cccc} 1 & 0 & 0 &0\\ 0& 0& 0 & 0\end{array}\right) \right\}$$

%AB=A3_5_a11a40b20_c_c11_d10d30_d40_d21 := subs([a_1 = 1, a_2 = 0, a_3 = 0, a_4 = 0, b_1 = 0, b_2 = 0, b_3 = 0, b_4 = 0, c_1 = 1, c_2 = 0, c_3 = 0, c_4 = 0, d_1 = 0, d_2 = 1, d_3 = 0, d_4 = 0], AX5)

$$D_{3,3}^3:=\left\{A=\left(\begin{array}{cccc} 1 & 0 & 0 &0\\ 0 & 0 & 0 &0\end{array}\right),\ \ B=\left(\begin{array}{cccc} 1 & 0 & 0 &0\\ 0& 1& 0 & 0\end{array}\right) \right\}.$$

The algebras $D_{3,3}^2$ and $D_{3,3}^3$ are not isomorphic because there is no an element of the automorphism group  $$Aut(As_{3,3}^2)=\left\{ \left(\begin{array}{crr} 1 & 0\\ z & t \end{array}\right)\right\}$$ sending $B$ of $D_{3,3}^2$ to $B$ of $D_{3,3}^3$.

The associative algebra $As_{3,3}^3=\begin{pmatrix}
 		1&0&0&0 \\
 		0& 1& 0& 0
 		\end{pmatrix}$ generates

%AB=A3_5_a11a40b21_c_c11_d21 := subs([a_1 = 1, a_2 = 0, a_3 = 0, a_4 = 0, b_1 = 0, b_2 = 1, b_3 = 0, b_4 = 0, c_1 = 1, c_2 = 0, c_3 = 0, c_4 = 0, d_1 = 0, d_2 = 1, d_3 = 0, d_4 = 0], AX5)

$$D_{3,3}^4:=\left\{A=\left(\begin{array}{cccc} 1 & 0 & 0 &0\\ 0 & 1 & 0 &0\end{array}\right),\ \ B=\left(\begin{array}{cccc} 1 & 0 & 0 &0\\ 0& 1& 0 & 0\end{array}\right) \right\}.$$
%AB=A3_5_a12_a40_b20_c_c12_d40_d30_d20 := subs([a_1 = 2, a_2 = 0, a_3 = 0, a_4 = 0, b_1 = 0, b_2 = 0, b_3 = 2, b_4 = 0, c_1 = 2, c_2 = 0, c_3 = 0, c_4 = 0, d_2 = 0, d_3 = 0, d_4 = 0], AX5)

%A3_a12_a40_b20 := subs([a_1 = 2, a_2 = 0, a_3 = 0, a_4 = 0, b_1 = 0, b_2 = 0, b_3 = 2, b_4 = 0], AX)
From $As_{3,3}^4=\begin{pmatrix}
 		2&0&0&0 \\
 		0&0& 2& 0
 		\end{pmatrix}$%A3_a12_a40_b22 := subs([a_1 = 2, a_2 = 0, a_3 = 0, b_1 = 0, b_2 = 2, b_3 = 2, b_4 = 0], AX) 
 \
 we get
$$D_{3,3}^5:=\left\{A=\left(\begin{array}{cccc} 2 & 0 & 0 &0\\ 0 & 0 & 2 &0\end{array}\right),\ \ B=\left(\begin{array}{cccc} 2 & 0 & 0 &0\\ \delta_1& 0& 0 & 0\end{array}\right), \ \delta_1 \in \mathbb{F} \right\}$$
%A3_5_a12_a40_b20_c_c12_d40_d30_d22d10 := subs([a_1 = 2, a_2 = 0, a_3 = 0, a_4 = 0, b_1 = 0, b_2 = 0, b_3 = 2, b_4 = 0, c_1 = 2, c_2 = 0, c_3 = 0, c_4 = 0, d_1 = 0, d_2 = 2, d_3 = 0, d_4 = 0], AX5)
$$D_{3,3}^6:=\left\{A=\left(\begin{array}{cccc} 2 & 0 & 0 &0\\ 0 & 0 & 2 &0\end{array}\right),\ \ B=\left(\begin{array}{cccc} 2 & 0 & 0 &0\\ 0& 2& 0 & 0\end{array}\right) \right\}$$
%A3_5_a12_a40_b20_c_c12_d40_d32_d := subs([a_1 = 2, a_2 = 0, a_3 = 0, a_4 = 0, b_1 = 0, b_2 = 0, b_3 = 2, b_4 = 0, c_1 = 2, c_2 = 0, c_3 = 0, c_4 = 0, d_3 = 0, d_1 = 0, d_2 = 0, d_3 = 2, d_4 = 0], AX5)
%A3_5_a12_a40_b20_c_c12_d40_d32_d := subs([a_1 = 2, a_2 = 0, a_3 = 0, a_4 = 0, b_1 = 0, b_2 = 0, b_3 = 2, b_4 = 0, c_1 = 2, c_2 = 0, c_3 = 0, c_4 = 0, d_1 = 0, d_2 = 0, d_3 = 2, d_4 = 0], AX5)
$$D_{3,3}^7:=\left\{A=\left(\begin{array}{cccc} 2 & 0 & 0 &0\\ 0 & 0 & 2 &0\end{array}\right),\ \ B=\left(\begin{array}{cccc} 2 & 0 & 0 &0\\ 0& 0& 2 & 0\end{array}\right) \right\}.$$
In order to check the isomorphisms between $D_{3,3}^5$, $D_{3,3}^6$ and $D_{3,3}^7$ we act by the elements of automorphism group
$$Aut(As_{3,3}^4)=\left\{ \left(\begin{array}{crr} 1 & 0\\ 1+2t & t \end{array}\right)\ t \neq 0 \right\}$$
to $B$ parts of each algebras:
$$\left(\begin{array}{cccc} 2 & 0 & 0 &0\\ \delta_1& 0& 0 & 0\end{array}\right)=g^{-1}\left(\begin{array}{cccc} 2 & 0 & 0 &0\\ 0& 2& 0 & 0\end{array}\right)g^{\otimes 2}, \ \mbox{where}\ g=\left(\begin{array}{lll} 1 & 0\\ 1+2t & t \end{array}\right) $$
$$\left(\begin{array}{cccc} 2 & 0 & 0 &0\\ 0& 2& 0 & 0\end{array}\right)=g^{-1}\left(\begin{array}{cccc} 2 & 0 & 0 &0\\ 0& 0& 2 & 0\end{array}\right)g^{\otimes 2}, \ \mbox{where}\ g=\left(\begin{array}{lll} 1 & 0\\ 1+2t & t \end{array}\right) $$
$$\left(\begin{array}{cccc} 2 & 0 & 0 &0\\ 0& 0& 2 & 0\end{array}\right)=g^{-1}\left(\begin{array}{cccc} 2 & 0 & 0 &0\\ \delta_1& 0& 0 & 0\end{array}\right)g^{\otimes 2}, \ \mbox{where}\ g=\left(\begin{array}{lll} 1 & 0\\ 1+2t & t \end{array}\right). $$
As a result we get inconsistent systems of equations.

Finally, $As_{3,3}^5=\begin{pmatrix}
 		2&0&0&\alpha_4 \\
 		0&2& 2& 0
 		\end{pmatrix}$ generates the following diassociative algebra
%A3_a12_a40_b22 := subs([a_1 = 2, a_2 = 0, a_3 = 0, b_1 = 0, b_2 = 2, b_3 = 2, b_4 = 0], AX)
%A3_5_a12_b22_c12c20c30c4a4d10d22d32d40 := subs([a_1 = 2, a_2 = 0, a_3 = 0, b_1 = 0, b_2 = 2, b_3 = 2, b_4 = 0, c_1 = 2, c_2 = 0, c_3 = 0, c_4 = a_4, d_1 = 0, d_2 = 2, d_3 = 2, d_4 = 0], AX5)
$$D_{3,3}^8:=\left\{A=\left(\begin{array}{cccc} 2 & 0 & 0 &\alpha_4\\ 0 & 2 & 2 &0\end{array}\right),\ \ B=\left(\begin{array}{cccc} 2 & 0 & 0 &\alpha_4\\ 0& 2& 2 & 0\end{array}\right) \right\}$$

\begin{theor} \label{char3Dias} Any non-trivial $2$-dimensional associative dialgebra over a field $\mathbb{F},$ $(Char(\mathbb{F})= 3)$ is isomorphic to only one of the following listed by their matrices of structure constants, such algebras:

\begin{enumerate}
%\item Diassociative algebras generated by $A_{12,2}^1$:
%a_1 = 0, a_2 = 0, a_3 = 0, a_4 = 0, b_1 = 1, b_2 = 0, b_3 = 0, b_4 = 0, c_1 = 0, c_2 = 0, c_3 = 0, c_4 = 0, d_2 = 0, d_3 = 0, d_4 = 0
%\begin{itemize}
\item $D_{13,3}^1:=\left\{A=\left(\begin{array}{cccc} 0 & 0 & 0 &0\\ 1 & 0 & 0 &0\end{array}\right),\ \ B=\left(\begin{array}{cccc} 0 & 0 & 0 &0\\ \delta_1& 0& 0 & 0\end{array}\right), \ \ \delta_1 \in \mathbb{F} \right\}$
\item $D_{3,3}^2:=\left\{A=\left(\begin{array}{cccc} 1 & 0 & 0 &0\\ 0 & 0 & 0 &0\end{array}\right),\ \ B=\left(\begin{array}{cccc} 1 & 0 & 0 &0\\ 0& 0& 0 & 0\end{array}\right) \right\}$

    %$$D_{6,2}^3:=\left\{A=\left(\begin{array}{cccc} 1 & 0 & 0 &0\\ 0 & 0 & 1 &0\end{array}\right)\ \ B=\left(\begin{array}{cccc} 1 & 0 & 0 &0\\ \beta_1& 0& 0 & 0\end{array}\right), \ \ \beta_1 \in \mathbb{F}\right\}$$
%AB=A6_5_a1n_a40_c_c1n1_d40_d30_d10d21 := subs([a_1 = -1, a_2 = 0, a_3 = 0, a_4 = 0, b_1 = 0, b_2 = 0, b_3 = -1, b_4 = 0, c_1 = -1, c_2 = 0, c_3 = 0, c_4 = 0, d_1 = 0, d_2 = 1, d_3 = 0, d_4 = 0], AX5)
%$$$$
%AB=A6_5_a1n_a40_c_c1n1_d40 := subs([a_1 = -1, a_2 = 0, a_3 = 0, a_4 = 0, b_1 = 0, b_2 = 0, b_3 = -1, b_4 = 0, c_1 = -1, c_2 = 0, c_3 = 0, c_4 = 0, d_1 = 0, d_2 = 0, d_3 = 1, d_4 = 0], AX5)
%$$$$
\item %a_1 = 1, a_2 = 0, a_3 = 0, a_4 = 0, b_1 = 0, b_2 = 0, b_3 = 0, b_4 = 0, c_1 = 1, c_2 = 0, c_3 = 0, c_4 = 0, d_1 = 0, d_2 = 1, d_3 = 0, d_4 = 0
$D_{3,3}^3:=\left\{A=\left(\begin{array}{cccc} 1 & 0 & 0 &0\\ 0 & 0 & 0 &0\end{array}\right),\ \ B=\left(\begin{array}{cccc} 1 & 0 & 0 &0\\ 0& 1& 0 & 0\end{array}\right) \right\}$
\item %a_1 = 1, a_2 = 0, a_3 = 0, a_4 = 0, b_1 = 0, b_2 = 1, b_3 = 0, b_4 = 0, c_1 = 1, c_2 = 0, c_3 = 0, c_4 = 0, d_1 = 0, d_2 = 1, d_3 = 0, d_4 = 0
$D_{3,3}^4:=\left\{A=\left(\begin{array}{cccc} 1 & 0 & 0 &0\\ 0 & 1 & 0 &0\end{array}\right),\ \ B=\left(\begin{array}{cccc} 1 & 0 & 0 &0\\ 0& 1& 0 & 0\end{array}\right) \right\}$

%$$$$
%A3_5_a12_a40_b20_c_c12_d40_d30_d22d10 := subs([a_1 = 2, a_2 = 0, a_3 = 0, a_4 = 0, b_1 = 0, b_2 = 0, b_3 = 2, b_4 = 0, c_1 = 2, c_2 = 0, c_3 = 0, c_4 = 0, d_1 = 0, d_2 = 2, d_3 = 0, d_4 = 0], AX5)
%$$$$
%A3_5_a12_a40_b20_c_c12_d40_d32_d := subs([a_1 = 2, a_2 = 0, a_3 = 0, a_4 = 0, b_1 = 0, b_2 = 0, b_3 = 2, b_4 = 0, c_1 = 2, c_2 = 0, c_3 = 0, c_4 = 0, d_3 = 0, d_1 = 0, d_2 = 0, d_3 = 2, d_4 = 0], AX5)
%A3_5_a12_a40_b20_c_c12_d40_d32_d := subs([a_1 = 2, a_2 = 0, a_3 = 0, a_4 = 0, b_1 = 0, b_2 = 0, b_3 = 2, b_4 = 0, c_1 = 2, c_2 = 0, c_3 = 0, c_4 = 0, d_1 = 0, d_2 = 0, d_3 = 2, d_4 = 0], AX5)
%$$$$
\item %a_1 = 1/2, a_2 = 0, a_3 = 0, a_4 = 0, b_1 = 0, b_2 = 0, b_3 = 1/2, b_4 = 0, c_1 = 1/2, c_2 = 0, c_3 = 0, c_4 = 0, d_2 = 0, d_3 = 0, d_4 = 0
$D_{3,3}^5:=\left\{A=\left(\begin{array}{cccc} 2 & 0 & 0 &0\\ 0 & 0 & 2 &0\end{array}\right),\ \ B=\left(\begin{array}{cccc} 2 & 0 & 0 &0\\ \delta_1& 0& 0 & 0\end{array}\right), \ \delta_1 \in \mathbb{F} \right\}$
\item %a_1 = 1/2, a_2 = 0, a_3 = 0, a_4 = 0, b_1 = 0, b_2 = 0, b_3 = 1/2, b_4 = 0, c_1 = 1/2, c_2 = 0, c_3 = 0, c_4 = 0, d_1 = 0, d_2 = 1/2, d_3 = 0, d_4 = 0
 $D_{3,3}^6:=\left\{A=\left(\begin{array}{cccc} 2 & 0 & 0 &0\\ 0 & 0 & 2 &0\end{array}\right),\ \ B=\left(\begin{array}{cccc} 2 & 0 & 0 &0\\ 0& 2& 0 & 0\end{array}\right) \right\}$

%$$$$

%AB=A3_5_a11_a40_b20_c_c11_d10d30_d40_d21 := subs([a_1 = 1, a_2 = 0, a_3 = 0, a_4 = 0, b_1 = 0, b_2 = 0, b_3 = 0, b_4 = 0, c_1 = 1, c_2 = 0, c_3 = 0, c_4 = 0, d_1 = 0, d_2 = 1, d_3 = 0, d_4 = 0], AX5)
%$$$$

\item %a_1 = 1/2, a_2 = 0, a_3 = 0, a_4 = 0, b_1 = 0, b_2 = 0, b_3 = 1/2, b_4 = 0, c_1 = 1/2, c_2 = 0, c_3 = 0, c_4 = 0, d_1 = 0, d_2 = 0, d_3 = 1/2, d_4 = 0
$D_{3,3}^7:=\left\{A=\left(\begin{array}{cccc} 2 & 0 & 0 &0\\ 0 & 0 & 2 &0\end{array}\right),\ \ B=\left(\begin{array}{cccc} 2 & 0 & 0 &0\\ 0& 0& 2 & 0\end{array}\right) \right\}$

\item %a_1 = 1/2, a_2 = 0, a_3 = 0, b_1 = 0, b_2 = 1/2, b_3 = 1/2, b_4 = 0, c_1 = 1/2, c_2 = 0, c_3 = 0, c_4 = a_4, d_1 = 0, d_2 = 1/2, d_3 = 1/2, d_4 = 0
$D_{3,3}^8:=\left\{A=\left(\begin{array}{cccc} 2 & 0 & 0 &\alpha_4\\ 0 & 2 & 2 &0\end{array}\right),\ \ B=\left(\begin{array}{cccc} 2 & 0 & 0 &\alpha_4\\ 0& 2& 2 & 0\end{array}\right) \right\}$
\end{enumerate}
\end{theor}
\section{Acknowledgement}
The author is gratefull to Professor U. Bekbaev for fruitifull disccussion on the research.


\begin{thebibliography}{9}
\bibitem{AdashevKhudoyberdievOmirov} {\sc  J.Q. Adashev, A.Kh. Khudoyberdiev, B.A. Omirov,}  Classification of complex naturally graded quasi-filiform Zinbiel algebras, \emph{Contemporary Math.}, AMS, 483, 2009.
\bibitem{AdashevOmirovKhudoyberdiev} {\sc  J.Q. Adashev, B.A. Omirov, A.Kh. Khudoyberdiev}, Classification of some classes of Zinbiel algebras, \emph{J. Gen. Lie Theor. and Appl.}, 4, 2010.
\bibitem{IJAC}	{\sc H. Ahmed, U.D. Bekbaev, I.S. Rakhimov}, Subalgebras, idempotents, ideals and quasi-units  of two-dimensional algebras, {\it International Journal of Algebra and Computations}, 2020, 30(5), 903--929. 
\bibitem{MJMS}	{\sc H. Ahmed, U.D. Bekbaev, I.S.Rakhimov}, Canonical form algorithm and separating systems for two-dimensional algebras, {\it Malaysian Journal of Math. Sciences}, 16(3): 413--435, 2022.
\bibitem{LJM}	{\sc H. Ahmed, U.D. Bekbaev, I.S.Rakhimov},  Identities on Two-Dimensional Algebras,
          {\it Lobachevskii Journal of Mathematics}, 2020, 41(9), 1615--1629.
\bibitem{AlbeverioAyupovOmirovKhudoyberdiyev} { \sc	 S. Albeverio, Sh.A. Ayupov, B.A. Omirov, A. Khudoyberdiyev}, $n$-dimensional filiform Leibniz algebras of length ($n$-1) and their derivations, \emph{Journal of Algebra}, 2008, 319, 2471--2488.
\bibitem{AlbeverioOmirovRakhimov2} { \sc S.	Albeverio, B.A. Omirov, I.S. Rakhimov,} Varieties of nilpotent complex Leibniz algebras of dimension less than five, \emph{Communications in Algebra}, 2005, 33, 1575--1585.

\bibitem {WI} {\sc W. Basri , I.M. Rikhsiboev}, On low-dimensional diassociative algebras, ICREM 3, UPM, 2007, 164--170.
\bibitem {RRB3} {\sc W. Basri,, I.S. Rakhimov,  I.M. Rikhsiboev}, On four-dimensional nilpotent Diassociative algebras, {\it Journal of Generalized Lie Theory and Applications}, 2015, 9(1).
\bibitem{Barnes2} { \sc Barnes D.W.}, On Levi's theorem for Leibniz algebras, \emph{Bull.
  Aust. Math. Soc.}, 86(2), 2012, 184--185.
\bibitem {UB} {\sc U. Bekbaev}, Classification of two-dimensional algebras over any basic field, 2023, {\it J. Phys.: Conf. Ser.}
\bibitem {UBm} {\sc U. Bekbaev}, On classification of $m$-dimensional algebras, 2017, {\it J. Phys.: Conf. Ser.}, 819 012012




\bibitem{CasasLadraOmirovKarimjanov2}{\sc  J.M. Casas, M. Ladra, B.A. Omirov, I.A. Karimjanov,}   Classification of solvable Leibniz algebras with naturally graded filiform nilradical, \emph{Linear Algebra and its Applications}, 2013, 438(7), 2973--3000.
\bibitem{CasasInsuaLadraLadra} {\sc J.M. Casas, M.A. Insua, M. Ladra, S. Ladra}, An algorithm for the classification of 3-dimensional complex Leibniz algebras,
\emph{Linear Algebra and its Applications}, 2012, 436(9), 3747--3756.
\bibitem{DjumadildaevTulenbaev} {\sc A.S. Dzhumadildaev, K.M. Tulenbaev,} Nilpotency of Zinbiel Algebras, \emph{Journal of Dynamical and Control Systems}, 2005, 11(2), 195--213.
\bibitem{GueIbrahimiFard} {\sc Guo L., Ebrahimi-Fard K.}, Rota-Baxter algebras and Dendriform algebras, \emph{Journal of Pure and Appl. Algebra}, 2008, 212, 320--339.
\bibitem{KhudRK}{\sc A.Kh. Khudoyberdiyev, I.S. Rakhimov, Sh.K. Said Husain},  On classification of 5-dimensional solvable Leibniz algebras, \emph{Linear Algebra and its Applications}, 2014, 457, 428--454.
\bibitem{LL} {\sc J.-L. Loday }, Une version non commutative des algebras de Lie: les algebras de Leibniz, \textit{Enseign. Math}, No. 39, (1993), 269-293.
\bibitem{LFCG}
{\sc J.-L. Loday, A. Frabetti, F.  Chapoton, F. Goichot},
\emph{Dialgebras and Related Operads},  Lecture Notes in
Mathematics, IV, 2001, 1763.
\bibitem{OR} {\sc	B.A. Omirov, I.S. Rakhimov}, On Lie-like filiform Leibniz algebras, {\it Bulletin of the Australian athematical Society}, 79: 391--404, (2009).
\bibitem{OR} {\sc	I.S. Rakhimov, U.D. Bekbaev}, On isomorphism classes and invariants of Finite Dimensional Complex Filiform Leibniz Algebras, {\it Communications in Algebra}, 38(12), 4705--4738.
\bibitem{RKhOM}{\sc I.S. Rakhimov, A.Kh.Khudoyberdiyev,  B.A.Omirov,  K.A.Mohd Atan}, On isomorphism criterion for a subclass of complex filiform Leibniz algebras, {\it International Journal of Algebra and Computations}, 2017, 27(7), 953-972.
\bibitem{R}{\sc I.S. Rakhimov}, On classification problem of Loday algebras, {\it Contemporary Mathematics}, AMS, 672, 225--244, 2016.
\bibitem {RRB1} {\sc I.M. Rikhsiboev, I.S. Rakhimov, W. Basri}, Classification of 3-dimensional complex diassociative algebras, {\it Malaysian Journal of Mathematical Sciences}, 2010, 4(2), 241--254.
\bibitem {RRB2} {\sc I.M. Rikhsiboev, I.S. Rakhimov,  W. Basri}, Diassociative algebras and their derivations, {\it Journal of Physics, Conference Series}, 2014, 553, 012006, 1--9.
\bibitem{IIW} {\sc Rikhsiboev I.M., Rakhimov I.S., Basri W.}, The description of dendriform algebra structures on two-dimensional complex space,  \emph{Journal of Algebra, Number Theory: Advances and Appl.}, 2010, 4(1), 1--18.
    \end{thebibliography}
\end{document}